\documentclass[12pt,a4paper]{article}
\usepackage[utf8]{inputenc}
\usepackage{amsmath, amssymb, amsthm, verbatim, hyperref}
\usepackage{mathrsfs}
\usepackage[dvipsnames]{xcolor}
\usepackage{ragged2e}
\usepackage{marginnote}
\usepackage{enumerate}
\usepackage{makecell}
\usepackage{longtable}
\usepackage{epsfig}
\usepackage{tikz}
\usepackage{ytableau}
\usepackage{hyperref}
\usepackage{caption}
\usepackage{subcaption}
\hypersetup{hidelinks}
\usepackage{fancyvrb}
\usepackage[OT1]{fontenc}
\usepackage{tikz-cd}
\usepackage[left=2.50cm, right=2.50cm, top=2.00cm, bottom=2.00cm]{geometry}
\usetikzlibrary {arrows.meta}
\usepackage{ifthen}

\newcommand{\ctab}[1]{\vcenter{\hbox{#1}}}

\def\Fl{\operatorname{Fl}}
\def\GL{\operatorname{GL}}
\def\FlVV{\Fl(V) \times \Fl(V)}

\def\tr{\operatorname{tr}}
\def\sh{\operatorname{sh}}

\def\evac{\operatorname{evac}}

\definecolor{goodgreen}{rgb}{0.01, 0.75, 0.24}

\theoremstyle{plain}

\newtheorem{innercustomthm}{Theorem}
\newenvironment{repthm}[1]
  {\renewcommand\theinnercustomthm{#1}\innercustomthm}
  {\endinnercustomthm}

\newtheorem{thm}{Theorem}[section]
\newtheorem{prop}[thm]{Proposition}
\newtheorem{cor}[thm]{Corollary}
\newtheorem{lemma}[thm]{Lemma}

\theoremstyle{definition}
\newtheorem{definition}[thm]{Definition}
\newtheorem{example}[thm]{Example}
\theoremstyle{remark}

\newtheorem{remark}[thm]{Remark}

\newcommand{\oRB}{\overline{RB}}

\newcommand{\SYT}{\operatorname{SYT}}
\newcommand{\OT}{\operatorname{OT}}

\newcommand{\gr}{\operatorname{gr}}
\newcommand{\im}{\operatorname{im}}
\newcommand{\bfrak}{\mathfrak{b}}

\newcommand{\cB}{\mathcal{B}}

\newcommand{\cF}{\mathcal{F}}
\newcommand{\cG}{\mathcal{G}}
\newcommand{\cH}{\mathcal{H}}
\newcommand{\cL}{\mathcal{L}}
\newcommand{\cN}{\mathcal{N}}
\newcommand{\cO}{\mathcal{O}}

\newcommand{\bfT}{\mathbf{T}}
\newcommand{\scrF}{\mathscr{F}}

\newcommand{\CC}{\mathbb{C}}

\newcommand{\BB}{\mathcal{B}}

\def\gfrak{\mathfrak{g}}
\newcommand{\eps}{\varepsilon}

\title{\vspace{-1em} 
\textbf{RSK with a twist, mirabolic Steinberg variety, and Pak--Postnikov representations}}
\author{\vspace{-2em} \\
Son Nguyen}
\date{\vspace{-3em}}

\allowdisplaybreaks

\begin{document}
\ytableausetup{centertableaux}

\maketitle

    \begin{abstract}
        There is a beautiful interplay between the classical Robinson--Schensted correspondence and the Steinberg variety, which is the union of conormal components of $\FlVV$. The irreducible components of the Steinberg variety gives the geometric interpretation of RS, and Pak--Postnikov description of RS as $U,D$ operators is the combinatorics shadow of the $S_n\times S_n$-representation inside the top homology of the Steinberg variety. We study the conormal components of $\FlVV\times L(\lambda)$, where $L(\lambda)$ is an irreducible representation of $\GL(V)$. We show that there are only three families of $L(\lambda)$ such that $\FlVV \times L(\lambda)$ has finitely many orbits. These families are determinant twists of classical RS, Travkin's mirabolic RSK and its Fourier dual. We extend the combinatorial interpretation of RS in terms of $U,D$ operators to mirabolic RSK. Finally, we give a filtration of the top homology of the mirabolic Steinberg variety in which each graded piece is a Pak--Postnikov $S_n \times S_n$-representations.
    \end{abstract}

\tableofcontents


\section{Introduction}

\subsection{The playbook}

    The \textit{Robinson--Schensted(--Knuth) correspondence}, introduced by Robinson \cite{robinson1938representations} and Schensted \cite{schensted1961longest}, is a classical operation in algebraic combinatorics. It is a simple bijection between \textit{permutations} and pairs of \textit{standard Young tableaux} of the same shape. It inspired several famous results in combinatorics, including plactic monoid \cite{fulton1997young}, which proves the Littlewood--Richardson rule, the Greene--Kleitman correspondence \cite{greene1976structure}, and several variations for other objects \cite{edelman1987balanced, tracy2004limit, haiman1989mixed,thomas2011longest,huang2022row,huang2022bumpless}.

    Behind the RS correspondence is a lesser-known but beautiful story, scattered all over the literature, that intertwines combinatorics and representation theory. In Section \ref{sec:RS}, we will document in details this playbook. The key ingredients are the following, and we refer the reader to Section \ref{sec:RS} for the notations.

    \begin{enumerate}
        \item On the representation theory side, start with the variety $X = \FlVV$, equipped with the diagonal $\GL(V)$ action. Its cotangent bundle has a canonical symplectic structure and a moment map $\mu: T^*X \rightarrow \gfrak^*$, where $\gfrak = \mathfrak{gl}(V)$. The zero fiber $Z = \mu^{-1}(0)$ is the classical Steinberg variety
        \[
        Z= \left\{(E,F,x)\in X \times\cN: xE_i \subseteq E_{i-1}, xF_i\subseteq F_{i-1}\right\}.
        \]
        \item Since $X$ has finitely many $\GL(V)$-orbits, the irreducible components of $Z$ are indexed by $\GL(V)$-orbits of $X$, which are in bijection with permutations $w\in S_n$. On the other hand, let $\lambda$ be the Jordan type of $x$, then the Springer fiber $\cB_x$ of $x$ has irreducible components indexed by SYTs of shape $\lambda$. Thus, one can also indexed the irreducible components of $Z$ by pairs of SYTs of the same shape. This is the geometric meaning of the RS correspondence.
        \item On the combinatorics side, one can interpret RS as $U$ and $D$ operators satisfying the simple relation $DU = UD + I$. This interpretation is deeply related to Fomin's growth diagram, which is interesting on its own and has inspired some beautiful combinatorial proofs.
        \item Pak and Postnikov used these operators to construct natural $S_p\times S_q$ representations $M(p,\beta,q)$ for all $p,q$ and sequence of nonzero integers $\beta$. When $p = q = n$ and $\beta = ()$, their representation is $\CC[S_n]$, which has the well-known Peter--Weyl decomposition
        \[
        \CC[S_n] = \bigoplus_{|\lambda| = n}\pi_\lambda \boxtimes \pi^*_\lambda,
        \]
        where $\pi_\lambda$ is the irreducible $S_n$-representation indexed by $\lambda$.
        \item Back to the representation theory side, the convolution product makes $H_{\operatorname{top}}(Z)$ an associative algebra that is isomorphic to $\CC[S_n]$. From this point of view, the Peter--Weyl decomposition witness the fact that if $x$ has Jordan type $\lambda$, then the top homology of the Springer fiber $H(\cB_x)$ is isomorphic to $\pi_\lambda$.
    \end{enumerate}

\subsection{Mirabolic RSK}

    Let $H_n$ be the Hecke algebra, which can be described as the Grothendieck ring of the constructible $\GL(V)$-equivariant mixed Tate complexes on $\FlVV$ with multiplication given by convolution. The two-sided Kazhdan-Lusztig cells are two-sided ideals spanned by the subsets of the Kazhdan-Lusztig basis formed by the classes of intersection cohomology sheaves, see \cite{kazhdan1979representations} and \cite{elias2020introduction} for an overview. Classical RS gives a (surprisingly) straightforward characterization of these two-sided cells: two $\GL(V)$-orbits on $\FlVV$ indexed by $w$ and $w'$ are in the same two-sided cell if and only if their RS have the same shape. These cells also appear in Lusztig’s classification of character sheaves on $\GL(V)$ \cite{lusztig2006character}, in which  each cell corresponds to a unique unipotent character sheaf.

    Now, let the Hecke algebra act by left and right convolution on the Grothendieck group of the constructible $\GL(V)$-equivariant mixed Tate complexes on $\FlVV\times V$. This comes with a Kazhdan-Lusztig basis indexed by the $\GL(V)$-orbits on $\FlVV\times V$, described in \cite{magyar1999multiple}. These orbits can be partitioned into bimodule KL cells similar to the classical Kazhdan--Lusztig cells. It is natural to ask if there is also a simple combinatorial characterization of these cells. Travkin's \textit{Mirabolic RSK} \cite{travkin2009mirabolic} is one such conjectural characterization. In particular, Mirabolic RSK is a bijection between $\GL(V)$-orbits on $\FlVV\times V$ and triples $(T_1,T_2,\theta)$, where $T_1,T_2$ are SYTs of shape $\nu,\nu'$ of size $n$, and $\nu/\theta,\nu'/\theta$ are horizontal strips. We denote by $\bfT_n$ be the set of such tripes $(\nu,\nu',\theta)$ where $|\nu| = |\nu'| = n$. Travkin conjectured that two orbits are in the same bimodule KL cell if and only if their mirabolic RSK have the same triple of shapes in $\bfT_n$ \cite[Conjecture 1]{travkin2009mirabolic}.

    Furthermore, analogous to Lusztig's classification of character sheaves on $\GL(V)$, in a sequel \cite{finkelberg2009mirabolic}, Finkelberg--Ginzburg--Travkin constructed explicit unipotent mirabolic character sheaves on $\GL(V) \times V$. They also conjectured that these  exhausts all the irreducible mirabolic character sheaves. This was later proved by Shoji--Sorlin in the finite-field setting \cite{shoji2014exotic}.

\subsection{Our results}

    One can apply this playbook for other varieties. In this paper, we apply it for $\FlVV \times L(\lambda)$, where $L(\lambda)$ is the irreducible representation of $\GL(V)$ indexed by $\lambda$.

    The first step is to identify for which $\lambda$, $\FlVV\times L(\lambda)$ has finitely many $\GL(V)$-orbits. We give a complete classification.

    \begin{repthm}{\ref{thm:classify}}
        Let $n\geq 2$ and let $L(\lambda)$ be an irreducible polynomial $\GL(V)$-representation. Then $\FlVV\times L(\lambda)$ has finitely many $G$-orbits if and only if $\lambda$ is one of the following
        \[
        \begin{array}{ll}
        (k^n),& k\geq1,\\
        (k+1,k^{n-1}),& k\geq0,\\
        ((k+1)^{n-1},k),& k\geq0.
        \end{array}
        \]
    \end{repthm}

    Let $D$ be the determinant representation and
    \[
    X_k = \Fl(V)^2 \times D^{\otimes k},\quad
    X^+_k = \Fl(V)^2 \times (V\otimes D^{\otimes k}),\quad
    X^-_k = \Fl(V)^2 \times (V^*\otimes D^{\otimes k+1})
    \]
    with moment maps $\mu_k,\mu^+_k$, $\mu^-_k$, respectively, and let
    \[
    Z_k = \mu_k^{-1}(0),\quad Z^+_k = (\mu^+_k)^{-1}(0),\quad Z^-_k = (\mu^-_k)^{-1}(0)
    \]
    be the zero fibers.
    
    The second step is to describe the geometric RSK correspondence. The correspondence for $Z_k$ is easy since it is the classical Steinberg variety with an extra vector $\delta \in D^{\otimes k}$.

    \begin{repthm}{\ref{thm:DkRSK}}
        There is a correspondence on the irreducible components of $Z_k$
        \[
        (w,\eps)\longleftrightarrow (P(w),Q(w),\eps),
        \]
        where $w\longleftrightarrow (P(w),Q(w))$ is the classical RSK, and $\eps = 0$ (resp. $\eps \neq 0$) if $\delta = 0$ (resp. $\delta \neq 0$).
    \end{repthm}

    For $Z^+_k$ and $Z^-_k$, it turns out that they are both isomorphic to the mirabolic Steinberg variety $Z_V$ in \cite{travkin2009mirabolic}. Thus, we get for free an RSK correspondence for $Z^+_k$ and $Z^-_k$ from Travkin's mirabolic RSK.

    \begin{repthm}{\ref{thm:Z+-}}
        We have $Z^+_k \simeq Z_V \simeq Z^-_k$, and hence their RSK is given by Travkin's mirabolic RSK.
    \end{repthm}

    Furthermore, for a vector space $W$, the canonical Fourier symplectic exchange is
    \[
    \cF_W:T^*W=W\times W^* \longrightarrow T^*(W^*)=W^*\times W,\qquad
    (w,\xi)\longmapsto(\xi,-w).
    \]
    This gives a Fourier bijection on the irreducible components of $T^*W$ and $T^*(W^*)$. It is the ``microlocal shadow'' of the sheaf-theoretic settings, where we have the Fourier--Sato transform for conic sheaves over $\CC$ \cite[Section 3.7]{kashiwara2013sheaves} and the Fourier--Deligne transform for $\ell$-adic sheaves \cite{laumon1987transformation, kiehl2001weil}. One  example is in Springer theory, where we have the Springer sheaf on the nilpotent cone and the Grothendieck sheaf on the semisimple locus, both of which carry the Weyl group action. The Fourier--Sato transform gives an isomorphism between these two sheaves \cite[Lemma 2.2]{achar2014weyl}. See also \cite{brylinski1982transformations, lusztig2010introduction, tamarkin2013microlocal, d2018microlocal}. For classical RS, Fourier exchange is the well-known involution $w\mapsto w_0ww_0$ that after RS sends $(P,Q)\mapsto (\evac(P),\evac(Q))$. For other instances, see \cite{travkin2009mirabolic, finkelberg2025lagrangian}.

    In our current context, we have two choices on which to perform the Fourier exchange, namely on $V$
    \[
    \cF:T^*V \longrightarrow T^*(V^*),\qquad
    (v,\alpha)\longmapsto(\alpha,-v),
    \]
    and on $\FlVV$
    \[
    \Delta : T^*(\FlVV) \longrightarrow T^*(\Fl(V^*) \times \Fl(V^*)), \qquad (E,F,x_1,x_2) \mapsto (E^\vee,F^\vee,-x_1^*,-x_2^*).
    \]
    Thus, we have a commuting diagram
    \[\begin{tikzcd}[sep = scriptsize]
    	{Z^+_V} && {Z^-_V} \\
    	\\
    	{Z^-_{V^*}} && {Z^+_{V^*}}
    	\arrow["{\cF}", from=1-1, to=1-3]
    	\arrow["\Delta"', from=1-1, to=3-1]
    	\arrow["\Delta", from=1-3, to=3-3]
    	\arrow["{\cF}"', from=3-1, to=3-3]
        \arrow["{\mathcal{T}}"{description}, dotted, from=1-1, to=3-3]
    \end{tikzcd}\]
    in which the diagonal is Travkin's Fourier exchange on the whole $\FlVV \times V$ \cite[Section 5.3]{travkin2009mirabolic}.

    The third step is to interpret mirabolic RSK as $U$ and $D$ operators. This is the topic of Section \ref{sec:mira-UD}.

    \begin{repthm}{\ref{thm:mira-UD}}
        Let $\pi$ be a partial permutation, let $(T_1,T_2,\theta)$ be the mirabolic RSK of $\pi$, and $(T_1',T_2',\theta')$ be image of $\pi$ under the $U$ and $D$ operators. Then,
        \[
        T_1 = \evac(T_1'),\qquad T_2 = \evac(T_2'),\qquad \theta = \theta'.
        \]
    \end{repthm}

    \begin{remark}
        In \cite{nguyen2026rsk}, we will give another growth diagram description for mirabolic RSK that does not require the additional translation to partial permutations and retains the geometric information similar to Fomin's growth diagram.
    \end{remark}

    This interpretation suggests that we can find Pak--Postnikov representations inside the top homology $H(Z_V)$. Let
    \[
    Z_{\leq k}=\{(E,F,v;x_1,x_2,\alpha)\in Z_V:\dim \CC[x_1]v\leq k\},\qquad Z_{\le-1}=\varnothing,
    \]
    and define the filtration
    \[
    \scrF_k H(Z_V)=\im \left(H(Z_{\leq k}) \longrightarrow H(Z_V)\right), \qquad -1\leq k\leq n.
    \]
    The graded pieces of this filtration are isomorphic to Pak--Postnikov representations.
    
    \begin{repthm}{\ref{thm:gr-PP}}
        For every $0\leq k\leq n$, there is an isomorphism of $S_n\times S_n$-representations
        \[
        \gr_k H(Z_V) \simeq M(n,(-k,k),n) \simeq \bigoplus_{(\nu,\nu',\theta)\in \bfT_n, |\theta| = n-k} \pi_\nu \boxtimes \pi^*_{\nu'}.
        \]
    \end{repthm}

    \begin{remark}
        It was pointed out to us by Jialiang Zou that the graded pieces in Theorem \ref{thm:gr-PP} are exactly the rank pieces of the type II oscillator bimodule in \cite[Section 4.3]{ma2025theta} when $m = n$ and $q = 1$. It will nice to explore the reason behind this under the frame of relative Langlands duality.
    \end{remark}
    
\section*{Acknowledgement}

    We thank Alex Postnikov for suggesting the Pak--Postnikov's $U$ and $D$ operators, which was the starting point of this project. We thank Jialiang Zou for introducing the mirabolic RSK to us. We also thank Daoji Huang and Linus Setiabrata for helpful conversations and remarks.

\section{Classical RS}\label{sec:RS}

\subsection{The moment map}

    First, we review the geometry behind RSK, aimed at a combinatorial audience. We refer the reader to \cite{da2004lectures} for a more complete introduction.

    Let $V = \CC^n$, $G = \GL(V)$ and $\gfrak = \mathfrak{gl}(V)$. Let $(M,\omega)$ be a holomorphic symplectic manifold with a holomorphic symplectic $G$-action. For $A \in \gfrak$, the \textit{fundamental vector field} on $M$ induced by $A$ is
    \[
    A_M:q\mapsto \left.\dfrac{d}{dt}\right|_{t=0} \exp(tA)\cdot q.
    \]
    This is the velocity of the curve $t\mapsto \exp(tA)\cdot q$ through $q$ at $t = 0$.

    \begin{definition}\label{def:moment-map}
        A smooth map $\mu:M\to\mathfrak g^*$ is a \textit{moment map} if it satisfies
        \begin{enumerate}
            \item[(a)] it is $G$-equivariant for the coadjoint action
            \[
             \mu(g\cdot m)=\operatorname{Ad}_g^*\mu(m),
            \]
            and
            \item[(b)] for every $A\in\mathfrak g$, the fundamental vector field $A_M$ satisfies
            \[
            d\mu^A=\omega(A_M,\cdot),
            \]
            where $\mu^A: m\mapsto \langle\mu(m),A\rangle$.
        \end{enumerate}
        The $G$-action is called \textit{Hamiltonian} if there exists such moment map.
    \end{definition}

    Let $X$ be a smooth complex variety with a $G$-action. The cotangent bundle
    \[
    M = T^*X = \{(q,\xi):q\in X,\ \xi\in T_q^*X\}
    \]
    has a canonical holomorphic symplectic structure. It is not hard (but also not easy) to check the following.

    \begin{lemma}\label{lem:mu}
        The moment map $\mu:T^*X \rightarrow \gfrak^*$ is given by
        \[
        \langle\mu(q,\xi),A\rangle=\xi(A_X(q)).
        \]
    \end{lemma}

    For example, we can check that $\mu$ is $G$-equivariant. The $G$-action on $T^*X$ is the cotangent lift given by
    \[
    g\cdot (q,\xi) = \left(g q,  \xi \circ (da_{g^{-1}})_{gq} \right),
    \]
    where $a_g:q\mapsto gq$. Thus,
    \begin{align*}
        \langle \mu(g\cdot(q,\xi)), A \rangle &= \langle \mu(g q,  \xi \circ (da_{g^{-1}})_{gq}), A \rangle \\
        &= \xi \circ (da_{g^{-1}})_{gq}(A_X(gq)) \\
        &= \xi ((\operatorname{Ad}_{g^{-1}}A)_X(q))\\
        &= \langle \mu(q,\xi), \operatorname{Ad}_{g^{-1}}A \rangle \\
        &= \langle \operatorname{Ad}^*_{g}\mu(q,\xi), A \rangle
    \end{align*}

    For a smooth locally closed subvariety $\mathcal O\subseteq X$, its conormal bundle is
    \[
    T^*_{\cO}X
    = \left\{(q,\xi):q\in\mathcal O,\ \xi|_{T_q\mathcal O}=0\right\}.
    \]
    We refer the readers to \cite[Section~3.3]{da2004lectures} for the Lagrangian property of conormal bundles. When $\cO$ is a $G$-orbit, $\xi|_{T_q(G\cdot q)}=0$ is equivalent to $\xi(A_X(q))=0$ for all $A\in \gfrak$. This means $\mu(q,\xi) = 0$. Hence,
    \[
    \mu^{-1}(0) = \bigcup_{\cO}\overline{T^*_{\cO}X}.
    \]
    If, in addition, $X$ is equidimensional and has finitely many $G$-orbits, then
    \[
    \bigl(\mu^{-1}(0)\bigr)_{\mathrm{red}}
    =\bigcup_{\mathcal O}\overline{T^*_{\mathcal O}X},
    \]
    and the irreducible components are exactly the closures indexed by $G$-orbits of $X$.

\subsection{RS correspondence}

    For the RS correspondence, let $X = \FlVV$. A \textit{(complete) flag} in $V$ is a sequence
    \[
    F = (0 = F_0 \subset F_1 \subset \cdots \subset F_n = V),
    \]
    where $\dim F_i = i$. We will use the following shorthand in our examples. Given a basis $e_1,\ldots,e_n$ of $V$, we say
    \[
    F = e_{i_1}~)~e_{i_2}~)~\cdots ~)~e_{i_n}
    \]
    if
    \[
    F = 0 \subset \langle e_{i_1} \rangle \subset \langle e_{i_1},e_{i_2} \rangle \subset \cdots \subset \langle e_{i_1},e_{i_2},\ldots, e_{i_n} \rangle.
    \]
    The Bruhat decomposition indexes $\GL(V)$-orbits on $X$ by permutations in $S_n$. In particular, a pair of flags $(E,F)$ is in the orbit $\cO_w$ if there is a basis $e_1,\ldots,e_n$ with
    \[
    E = e_{1}~)~e_{2}~)~\cdots ~)~e_{n},\qquad
    F = e_{w(1)}~)~ e_{w(2)}~)~\cdots~)~e_{w(n)}.
    \]
    Let $\cN\subset \gfrak$ be the nilpotent cone. We identify $\gfrak$ with $\gfrak^*$ by the trace pairing $\langle x,A\rangle=\tr(xA)$. Then, the cotangent bundle of $X$ is
    \[
    T^*X = \left\{(E,F,x_1,x_2)\in X \times\cN \times\cN: x_1E_i \subseteq E_{i-1}, x_2F_i\subseteq F_{i-1}\right\}.
    \]
    For a flag $F$, let
    \[
    \bfrak_F=\{A\in\gfrak:AF_i\subseteq F_i\text{ for all }i\}
    \]
    be its Borel subalgebra. Note that if $xF_i\subset F_{i-1}$ for all $i$, then $\tr (x\bfrak_F) = 0$. At $q = (E,F)$, we have $A_X(q) = (A + \bfrak_E,A + \bfrak_F)$. Hence, the moment map is given by
    \[
    \langle\mu(E,F,x_1,x_2),A\rangle=(x_1,x_2)(A + \bfrak_E, A + \bfrak_F)
    \]
    \[
    = \tr x_1(A + \bfrak_E) + \tr x_2 (A + \bfrak_F) = \langle x_1 + x_2, A \rangle.
    \]
    Thus, $\mu(E,F,x_1,x_2) = x_1+x_2$. Let $Z = \mu^{-1}(0)$, then
    \[
    Z = \left\{(E,F,x)\in X \times\cN: xE_i \subseteq E_{i-1}, xF_i\subseteq F_{i-1}\right\}.
    \]
    This is the \textit{Steinberg variety}, and the irreducible components of $Z$ are indexed by $w\in S_n$.

    On the other hand, let $x\in \cN$ be a nilpotent operator of Jordan type $\nu$. The \textit{Springer fiber} of $x$
    \[
     \BB_x(V)=\{F\in \Fl(V):xF_i\subseteq F_{i-1}\text{ for all }i\}
    \]
    has irreducible components indexed by the set $\SYT(\nu)$ of \textit{standard Young tableaux (SYT)} of shape $\nu$. Let $F$ be a generic flag in the component indexed by a SYT $T$, the Jordan types $J(x|_{F_i})$ are the shapes of the subtableaux with entries $1,\ldots,i$. Hence, one can also index the irreducible components of $Z$ by the irreducible components of the two Springer fibers, which is a pair of SYTs of the same shape.

    As a result, the irreducible components of $Z$ give a bijection between permutations $w\in S_n$ and pairs of SYTs $(P,Q)$ of the same shape and size $n$. Steinberg \cite{steinberg1988occurrence} showed that this is the RS correspondence, which can be describe combinatorially as follows. We also refer the readers to \cite{stanley2023enumerative} for a more comprehensive combinatorial introduction.

    \begin{definition}
        Let $w=w(1)\cdots w(n)\in S_n$. The
        \textit{Robinson--Schensted correspondence} assigns to $w$ a pair $(P(w),Q(w))$ of SYTs of the same shape, called the \textit{insertion tableau} and the \textit{recording tableau}, as follows.
        
        Start with $P=Q=\varnothing$. For $i=1,\ldots,n$, insert $w(i)$
        into $P$ by \emph{row insertion}.
        \begin{enumerate}
            \item If $w(i)$ is larger than every entry in row 1, append $w(i)$ to the end of row 1.
            \item Otherwise, in row 1, replace the leftmost entry larger than $w(i)$ by $w(i)$, and insert the replaced entry
            into the next row by the same rule.
            \item This creates exactly one new box in $P$. Add a box in the same position in $Q$ and fill it with $i$.
        \end{enumerate}
        After all $n$ insertions, set $P(w)=P$ and $Q(w)=Q$.
    \end{definition}

    \begin{example}\label{ex:RS}
        Let $w = 2341$, we build the insertion tableau as follows
        \[
        \ctab{\begin{ytableau} 
            *(yellow) 2
        \end{ytableau}}
        ~\longrightarrow~
        \ctab{\begin{ytableau} 
            2 & *(yellow) 3
        \end{ytableau}}
        ~\longrightarrow~
        \ctab{\begin{ytableau} 
            2 & 3 & *(yellow) 4
        \end{ytableau}}
        ~\longrightarrow~
        \ctab{\begin{ytableau} 
            1 & 3 & 4 \\
            *(yellow) 2
        \end{ytableau}} ~=~ P.
        \]
        The yellow boxes are the new boxes in each step, so
        \[
        Q ~=~ \ctab{\begin{ytableau} 
            1 & 2 & 3 \\
            4
        \end{ytableau}}.
        \]
        Let us check the geometry. Let
        \[
        E = e_{1}~)~e_{2}~)~e_3~)~e_{4},\qquad
        F = e_{2}~)~ e_{3}~)~e_4~)~e_1.
        \]
        A generic nilpotent operator $x$ that stabilizes both $E$ and $F$ has the form
        \[
        \left( \begin{matrix}
            0 & 0 & 0 & 0 \\
            0 & 0 & E_{23} & E_{24} \\
            0 & 0 & 0 & E_{34} \\
            0 & 0 & 0 & 0
        \end{matrix} \right).
        \]
        For example, one can replace the nonzero entries by prime numbers. The sequence of Jordan types on $E$ is
        \[
        \ctab{\begin{ytableau} 
            {}
        \end{ytableau}}
        ~\subset~
        \ctab{\begin{ytableau} 
            {} \\
            {}
        \end{ytableau}}
        ~\subset~
        \ctab{\begin{ytableau} 
            {} & {} \\
            {}
        \end{ytableau}}
        ~\subset~
        \ctab{\begin{ytableau} 
            {} & {} & {} \\
            {}
        \end{ytableau}}.
        \]
        One can check that this is $P$. On $F$, the sequence is
        \[
        \ctab{\begin{ytableau} 
            {}
        \end{ytableau}}
        ~\subset~
        \ctab{\begin{ytableau} 
            {} & {}
        \end{ytableau}}
        ~\subset~
        \ctab{\begin{ytableau} 
            {} & {} & {}
        \end{ytableau}}
        ~\subset~
        \ctab{\begin{ytableau} 
            {} & {} & {} \\
            {}
        \end{ytableau}},
        \]
        which is $Q$.
    \end{example}

\subsection{Growth diagram}

    Fomin's growth diagram \cite{fomin1988generalized} is another description of the RS correspondence. See also \cite{britz1999finite,van2000flag}.

    \begin{definition}\label{def:GD}
        For $w\in S_n$, take an $n\times n$ grid with vertex coordinates
        $(i,j)$, where $i$ increases downwards and $j$ to the right.
        Mark square $(i,j)$ by $\times$ when $w(j)=i$, and put
        $\varnothing$ at every vertex on the North and West boundaries.
        For one square
        \[\begin{tikzcd}[sep = small]
        	{\lambda} && {\mu} \\
        	\\
        	{\nu} && {\rho}
            \arrow[no head, from=1-1, to=1-3]
        	\arrow[no head, from=1-1, to=3-1]
        	\arrow[no head, from=1-3, to=3-3]
        	\arrow[no head, from=3-1, to=3-3]
        \end{tikzcd},\]
        one determine $\rho$ from $\lambda,\mu,\nu$ by
        \[
        \rho=
        \begin{cases}
        \lambda+b_1,&\text{the square is marked},\\
        \mu\cup\nu,&\text{the square is unmarked and }\mu\ne\nu,\\
        \lambda,&\text{the square is unmarked and }\lambda=\mu=\nu,\\
        \mu+b_{k+1},&\text{the square is unmarked and }\mu=\nu=\lambda+b_k.
        \end{cases}
        \]
        Here, $b_k$ denotes a box in row $k$, and the union is taken row by row. Then, East and South boundary sequences give $P(w)$ and $Q(w)$,
        respectively.
    \end{definition}

    The growth diagram description clearly shows the symmetry of RS: if $w \leftrightarrow (P,Q)$ then $w^{-1} \leftrightarrow (Q,P)$. In addition, each vertex in the growth diagram retains the geometry information. Let $\lambda_{i,j}$ be the shape at vertex $(i,j)$. Then, by the previous section,
    \[
    \lambda_{i,n} = J(x|_{E_i}), \qquad \lambda_{n,i} = J(x|_{F_i}).
    \]
    This is actually true for all vertices, that is
    \[
    \lambda_{i,j} = J(x|_{E_i\cap F_j}).
    \]

    \begin{example}
        Let $w= 2341$ as in Example \ref{ex:RS}. The corresponding growth diagram is
        \begin{center}
            \includegraphics[scale = 1]{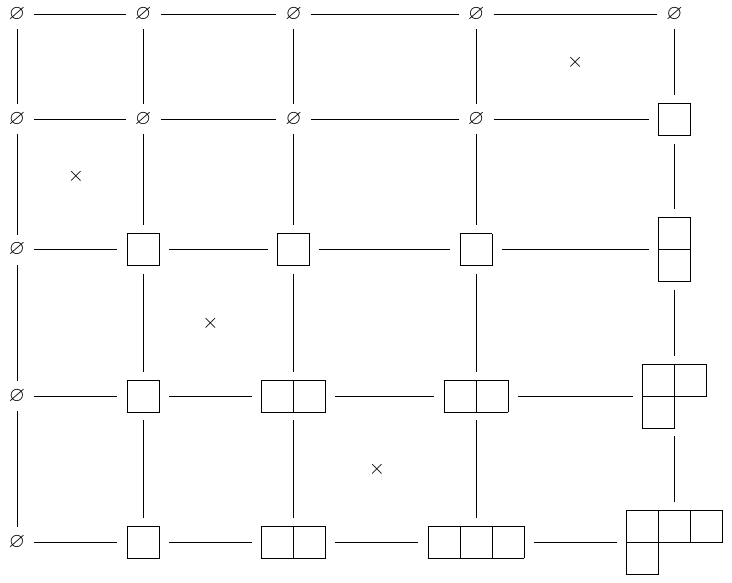}.
        \end{center}
        The shape $\lambda_{3,2} = (2)$, and one can check that
        \[
        J(x|_{E_3\cap F_2}) = J(x|_{\langle e_2,e_3\rangle})) = (2).
        \]
    \end{example}

\subsection{Up and down operators and Pak--Postnikov representations}\label{subsec:PPUD}

    In \cite{pak1996oscillating}, Pak--Postnikov observed that all four rules in Definition \ref{def:GD} can be written as a single rule in terms of toggle.

    \begin{definition}
        A (column-strict) \textit{oscillating tableau} of shape $(\lambda,\mu)$ is a sequence of partitions
        \[
        T = (\lambda = \lambda_0, \lambda_1,\lambda_2,\ldots, \lambda_\ell = \mu)
        \]
        such that for all $i$, either $\lambda_{i+1}/\lambda_i$ or $\lambda_i/\lambda_{i+1}$ is a (nonempty) horizontal strip. We represent $T$ as follows: if $\lambda_{i-1}\subset \lambda_{i}$, then we put an $i$ for each box added; if $\lambda_{i}\subset \lambda_{i-1}$, then we put an $-i$ for each box removed. The \textit{content} of $T$ is $\beta = (\beta_1,\ldots,\beta_\ell)$, where either $\lambda_i / \lambda_{i-1}$ is a horizontal $\beta_i$-strip or $\lambda_{i-1}/\lambda_i$ is a horizontal $(-\beta_i)$-strip.

        Let $\OT(\lambda,\mu,\beta)$ be the set of oscillating tableau of shape $(\lambda,\mu)$ and content $\beta$.
    \end{definition}

    \begin{example}
        The following oscillating tableau
        \ytableausetup{boxsize=2.5em}
        \[
        T = ~\ctab{\begin{ytableau}
            \none & \none & 1 & {1\overline{2}3\overline{5}} \\
            \none & {1\overline{2}4} & 4 \\
            1 \\
            3
        \end{ytableau}}
        \]
        is an oscillating tableau of shape $((2,1),(3,3,1,1))$ and content $(4,-2,2,2,-1)$.
    \end{example}

    \begin{example}
        Each pair of SYTs $(P,Q)$ of the same shape can be viewed as an oscillating tableau of shape $(\varnothing,\varnothing)$ and content $(1^n,(-1)^n)$ by adding the boxes in $P$ and removing the boxes in $Q$. For example, the pair in Example \ref{ex:RS} can be view as
        \ytableausetup{boxsize=2em}
        \[
        (P,Q) = ~\ctab{\begin{ytableau}
            1\overline{8} & 3\overline{7} & 4\overline{6} \\
            2\overline{5}
        \end{ytableau}}.
        \]
    \end{example}

    An oscillating tableau can be viewed as a walk in the Young's lattice, in which each step adds or removes a horizontal strip. Thus, we define the following operators.

    \begin{definition}
        Let $R$ be the infinite-dimensional space of formal linear combinations of partitions. Let
        \[
        U_n:\lambda \mapsto \sum_{\text{$\mu/\lambda$ is a horizontal $n$-strip}}\mu,\qquad D_n:\lambda \mapsto \sum_{\text{$\lambda/\mu$ is a horizontal $n$-strip}}\mu.
        \]
        We denote $U = U_1$ and $D = D_1$ be the operators that walk up and down the edges of the Young's lattice.
    \end{definition}

    Let $\langle b \rangle$ be $U_b$ if $b>0$ and $D_b$ if $b < 0$. Then, clearly
    \[
    |\OT(\lambda,\mu,\beta)| = [\mu](\langle b_\ell\rangle\cdots \langle b_1\rangle \lambda).
    \]
    The up and down operators have the following simple relations.

    \begin{thm}[{\cite[Theorem 6.1]{pak1996oscillating}}]\label{thm:UDrel}
        We have
        \[
        [U_m,U_n] = [D_m,D_n] = 0,\qquad [D_m,U_n] = D_{m-1}U_{n-1}.
        \]
    \end{thm}

    The first relation $[U_m,U_n] = [D_m,D_n] = 0$ is the Bender--Knuth involutions. Our main concern is the second relation, which can be proved by the following bijection.

    Define $\Phi_{m,n}:\OT(\lambda,\mu, (m,-n)) \longrightarrow \bigsqcup_{0\leq k \leq \min(m,n)}\OT(\lambda,\mu, (-n+k, m-k))$ as follows. An oscillating tableau $T$ in $\OT(\lambda,\mu, (m,-n))$ is a sequence $(\lambda,\nu,\mu)$ where $\nu/\lambda$ and $\nu/\mu$ are horizontal strips. This means that
    \[\begin{tikzcd}[sep = scriptsize]
    	& {\lambda_1} && {\lambda_2} && {\lambda_3} && \cdots \\
    	{\nu_1} && {\nu_2} && {\nu_3} && \cdots \\
    	& {\mu_1} && {\mu_2} && {\mu_3} && \cdots
    	\arrow["\geq"{description}, from=1-2, to=2-3]
        \arrow["\geq"{description}, from=1-4, to=2-5]
        \arrow["\geq"{description}, from=1-6, to=2-7]
        \arrow["\geq"{description}, from=2-1, to=1-2]
        \arrow["\geq"{description}, from=2-1, to=3-2]
        \arrow["\geq"{description}, from=2-3, to=1-4]
        \arrow["\geq"{description}, from=2-3, to=3-4]
        \arrow["\geq"{description}, from=2-5, to=1-6]
        \arrow["\geq"{description}, from=2-5, to=3-6]
        \arrow["\geq"{description}, from=2-7, to=1-8]
        \arrow["\geq"{description}, from=2-7, to=3-8]
        \arrow["\geq"{description}, from=3-2, to=2-3]
        \arrow["\geq"{description}, from=3-4, to=2-5]
        \arrow["\geq"{description}, from=3-6, to=2-7]
    \end{tikzcd}.\]
    We can toggle as follows: define $\nu'$ by $\nu'_i = \min(\lambda_i,\mu_i) + \max(\lambda_{i+1},\mu_{i+1}) - \nu_{i+1}$. One can check that this new sequence $(\lambda,\nu',\mu)$ is in $\OT(\lambda,\mu (-n+k, m-k))$ for $k = \nu_1 - \max(\lambda_1,\mu_1)$.

    \begin{example}
        Let
        \[
        T = ~\ctab{\begin{ytableau}
            \none & \none & \none & \none & 1 & 1\overline{2} \\
            \none & \none & \none & 1 \\
            \none & \none &  1 \\
            \overline{2} & 1\overline{2}
        \end{ytableau}} ~\in \OT((4,3,2,1),(5,4,3),(5,-3)).
        \]
        Here, $\nu = (6,4,3,2)$, and we have
        \[\begin{tikzcd}[sep=scriptsize]
        	& 4 && 3 && 2 && 1 \\
        	6 && 4 && 3 && 2 \\
        	& 5 && 4 && 3 && 0
        	\arrow["\geq"{description}, from=1-2, to=2-3]
        	\arrow["\geq"{description}, from=1-4, to=2-5]
        	\arrow["\geq"{description}, from=1-6, to=2-7]
        	\arrow["\geq"{description}, from=2-1, to=1-2]
        	\arrow["\geq"{description}, from=2-1, to=3-2]
        	\arrow["\geq"{description}, from=2-3, to=1-4]
        	\arrow["\geq"{description}, from=2-3, to=3-4]
        	\arrow["\geq"{description}, from=2-5, to=1-6]
        	\arrow["\geq"{description}, from=2-5, to=3-6]
        	\arrow["\geq"{description}, from=2-7, to=1-8]
        	\arrow["\geq"{description}, from=2-7, to=3-8]
        	\arrow["\geq"{description}, from=3-2, to=2-3]
        	\arrow["\geq"{description}, from=3-4, to=2-5]
        	\arrow["\geq"{description}, from=3-6, to=2-7]
        \end{tikzcd}.
        \]
        Toggling $\nu$ gives
        \[\begin{tikzcd}[sep=scriptsize]
        	4 && 3 && 2 && 1 & \\
        	& 4 && 3 && 1 && 0 \\
        	5 && 4 && 3 && 0
        	\arrow["\geq"{description}, from=1-1, to=2-2]
        	\arrow["\geq"{description}, from=1-3, to=2-4]
        	\arrow["\geq"{description}, from=1-5, to=2-6]
        	\arrow["\geq"{description}, from=1-7, to=2-8]
        	\arrow["\geq"{description}, from=2-2, to=1-3]
        	\arrow["\geq"{description}, from=2-2, to=3-3]
        	\arrow["\geq"{description}, from=2-4, to=1-5]
        	\arrow["\geq"{description}, from=2-4, to=3-5]
        	\arrow["\geq"{description}, from=2-6, to=1-7]
        	\arrow["\geq"{description}, from=2-6, to=3-7]
        	\arrow["\geq"{description}, from=3-1, to=2-2]
        	\arrow["\geq"{description}, from=3-3, to=2-4]
        	\arrow["\geq"{description}, from=3-5, to=2-6]
        	\arrow["\geq"{description}, from=3-7, to=2-8]
        \end{tikzcd},\]
        which is the oscillating tableau
        \[
        T' = ~\ctab{\begin{ytableau}
            \none & \none & \none & \none & 2 \\
            \none & \none & \none & 2 \\
            \none & \overline{1}2 &  2 \\
            \overline{1}
        \end{ytableau}} ~\in \OT((4,3,2,1),(5,4,3),(-2,4)).
        \]
        Here, $k = 6 - \max(4,5) = 1$.
    \end{example}

    In Theorem \ref{thm:UDrel}, when $m = 1$ or $n = 1$, we have
    \[
    DU_n = U_nD + U_{n-1}, \qquad D_mU = UD_m + D_{m-1}.
    \]
    When both $m = n = 1$, we have the well-known relation
    \[
    DU = UD + I.
    \]
    From this relation, one actually have $D^nU^n \varnothing = n!\varnothing$. This is the RS correspondence.
    
    Indeed, one can check that the cases for $\Phi_{1,1}$ are exactly the four cases in Definition \ref{def:GD}. Hence, the bijection $\Phi_{1,1}$ gives the inverse RS correspondence via the growth diagram as follows. Given $(P,Q)\in \OT(\varnothing,\varnothing,(1^n,(-1)^n))$, where $P$ records the $U$-steps and $Q$ records the $D$-steps. Put $P$ and $Q$ on the East and South boundary of the growth diagram. For each square
    \[\begin{tikzcd}[sep = small]
        {\lambda} && {\mu} \\
        \\
        {\nu} && {\rho}
        \arrow[no head, from=1-1, to=1-3]
        \arrow[no head, from=1-1, to=3-1]
        \arrow[no head, from=1-3, to=3-3]
        \arrow[no head, from=3-1, to=3-3]
    \end{tikzcd},\]
    we have $(\mu,\rho,\nu)\in \OT(\mu,\nu,(1,-1))$. Then let $\lambda$ be such that $\Phi_{1,1}(\mu,\rho,\nu) = (\mu,\lambda,\nu)$. If $(\mu,\lambda,\nu) \in \OT(\mu,\nu,())$, then put an $\times$ in the square. Then, the $\times$'s give a permutation in $S_n$.

    In addition to the up and down operators, Pak--Postnikov also introduced a natural $S_p\times S_q$ representation $M(p,\beta,q)$ on the space of intransitive graphs. We refer the readers to \cite[Section 3]{pak1996oscillating} for the construction. This representation has the following important decomposition.

    \begin{thm}[{\cite[Theorem 3.1]{pak1996oscillating}}]\label{thm:PPrep}
        Let $\pi_\lambda$, where $|\lambda| = n$ be the irreducible representation of $S_n$ indexed by $\lambda$. We have
        \[
        M(p,\beta,q) = \bigoplus_{|\lambda| = p,|\mu| = q} |\OT(\lambda,\mu,\beta)| \pi_\lambda \boxtimes \pi^*_\mu.
        \]
    \end{thm}

    \begin{example}\label{ex:PPrep-pp}
        An important special case of Theorem \ref{thm:PPrep} is when $p = q = n$ and $\beta = ()$. Then, $M(n,(),n) \simeq \CC[S_n]$, so Theorem \ref{thm:PPrep} gives the well-known Peter--Weyl decomposition
        \[
        \CC[S_n] = \bigoplus_{|\lambda| = n}\pi_\lambda \boxtimes \pi^*_\lambda.
        \]
    \end{example}

\subsection{Homology of Steinberg variety}\label{subsec:convo-hom-st}

    Let $x$ be a nilpotent operator of Jordan type $\lambda$, then $H(\cB_x) = \pi_\lambda$. Thus, it is natural to ask whether one can see the Peter--Weyl decomposition from the Steinberg variety point of view. This is also a classical story that we refer the readers to \cite[Section 3]{chriss1997representation} for a comprehensive introduction.

    Let $V = \CC^n$, recall that the Steinberg variety is
    \[
    Z = \left\{(E,F,x)\in \FlVV \times\cN: xE_i \subseteq E_{i-1}, xF_i\subseteq F_{i-1}\right\}.
    \]
    Let $H(Z)$ be the top Borel--Moore homology. Then,
    \[
    H(Z) = \bigoplus_w \CC\cdot[Z_w],
    \]
    where $Z_w = \overline{T^*_{\cO_w}\FlVV}$. Next, we need to define the convolution product on $H_*(Z)$. Let $M_1,M_2,M_3$ be connected, oriented $C^\infty$-manifolds and
    \[
    Z_{12} \subset M_1\times M_2,\qquad Z_{23} \subset M_2\times M_3
    \]
    be closed subsets. Define the composition
    \[
    Z_{12}\circ Z_{23} = \{(m_1,m_3) \in M_1\times M_3~|~\exists~ m_2\in M_2:(m_1,m_2) \in Z_{12}, (m_2,m_3) \in Z_{23}\}.
    \]
    Let $d = \dim M_2$, then we can define a convolution in Borel--Moore homology by
    \[
    H_i(Z_{12})\times H_j(Z_{23}) \longrightarrow H_{i+j-d}(Z_{12}\circ Z_{23}), \qquad (c_{12},c_{23}) \mapsto c_{12}\ast c_{23}.
    \]
    Now let $M$ be a smooth complex manifold and let $M_1 = M_2 = M_3 = M$. Let $N$ be a variety and a proper map $s:M \rightarrow N$. Let
    \[
    Z_{12} = Z_{23} = Z = M\times_N M = \{(m_1,m_2) \in M\times M~|~s(m_1) = s(m_2)\}.
    \]
    Then $Z\circ Z = Z$, so we have a convolution map
    \[
    H_*(Z)\times H_*(Z) \longrightarrow H_*(Z).
    \]

    \begin{cor}
        $H_*(Z)$ has a natural structure of an associative algebra with unit given by the fundamental class of $[M_{\Delta}]$. In addition, the top homology $H(Z)$ is a subalgebra of $H_*(Z)$.
    \end{cor}
    
    Now let us get back to our setting. Let $M = \tilde{\cN}$, $N = \cN$, and $s:\tilde{\cN} \rightarrow \cN$ be the Springer resolution. Then, the convolution product gives $H_*(Z)$ a natural structure of an associative algebra in which $H(Z)$ is a subalgebra. In fact, we have more: there is an algebra isomorphism
    \[
    H(Z) \simeq \CC[S_n] = \bigoplus_{|\lambda| = n}\pi_\lambda \boxtimes \pi^*_\lambda.
    \]

\section{Mirabolic RSK}\label{subsec:miraRSK}

    \ytableausetup{boxsize=1.5em}

    Now, let $X = \FlVV \times V$. Travkin \cite[Section 3.1]{travkin2009mirabolic} showed that the $\GL(V)$-orbits on $X$ can be indexed by $RB$ permutations.
    
    \begin{definition}\label{def:RB}
        An $RB$ permutation is a pair $\Tilde{w} = (w,\beta)$, where $w\in S_n$ and $\beta \subset [n]$ such that if $i\notin \beta$ and $j \in \beta$, then either $i > j$ or $w(i) > w(j)$. Let $RB_n$ be the set of RB permutations with $w\in S_n$.
    \end{definition}

    \begin{example}\label{ex:RB-perm}
        We have $(2431,\{1,4\}) \in RB_4$. One can draw $\Tilde{w}$ in an $n\times n$ grid by putting an $\times$ in $(w(i),i)$ and circle the $\times$'s in all columns in $\beta$. For example, $(2431,\{1,4\})$ is
        \[
        \begin{ytableau} 
            {} & {} & {} & {\otimes} \\
            {\otimes} & {} & {} & {} \\
            {} & {} & {\times} & {} \\
            {} & {\times} & {} & {}
        \end{ytableau}.
        \]
        Note that the condition ``if $i\notin \beta$ and $j \in \beta$ then either $i > j$ or $w(i) > w(j)$'' means that the $\otimes$'s are NW closed. For instance, $(2431,\{2,4\})$ is
        \[
        \begin{ytableau} 
            {} & {} & {} & {\otimes} \\
            {\times} & {} & {} & {} \\
            {} & {} & {\times} & {} \\
            {} & {\otimes} & {} & {}
        \end{ytableau},
        \]
        so this pair is not an RB permutation since the $\times$ in column $1$ is NW of the $\otimes$ in column $2$.
    \end{example}

    A triple $(E,F,v) \in \FlVV \times V$ is in the orbit $\cO_{\Tilde{w}}$ if there is a basis $e_1,\ldots,e_n$ with
    \[
    E = e_{1}~)~e_{2}~)~\cdots ~)~e_{n},\qquad
    F = e_{w(1)}~)~ e_{w(2)}~)~\cdots~)~e_{w(n)}, \qquad v= \sum_{i\in\beta} e_{w(i)}.
    \]

    \begin{example}\label{ex:RB-orbit}
        Continuing Example \ref{ex:RB-perm}, the following triple is in $\cO_{(2431,\{1,4\})}$
        \[
        E = e_{1}~)~e_{2}~)~e_3 ~)~e_4,\qquad
        F = e_{2}~)~ e_{4}~)~e_3~)~e_{1}, \qquad v= e_1 + e_2.
        \]
    \end{example}
    
    The cotangent bundle of $X$ is
    \[
    T^*X = \left\{(E,F,v,x_1,x_2,\alpha)\in X \times\cN \times\cN \times V^*: x_1E_i \subseteq E_{i-1}, x_2F_i\subseteq F_{i-1}\right\}.
    \]
    At $q = (E,F,v)$, we have $A_X(q) = (A + \bfrak_E,A + \bfrak_F, Av)$. Hence, the moment map is given by
    \[
    \langle\mu(E,F,v,x_1,x_2,\alpha),A\rangle=(x_1,x_2,\alpha)(A + \bfrak_E, A + \bfrak_F,Av)
    \]
    \[
    = \tr x_1(A + \bfrak_E) + \tr x_2 (A + \bfrak_F) + \tr (v\otimes \alpha)A = \langle x_1 + x_2 + v\otimes \alpha, A \rangle.
    \]
    Thus, $\mu(E,F,v,x_1,x_2,\alpha) = x_1+x_2 + v\otimes \alpha$. The \textit{mirabolic Steinberg variety} is
    \[
    Z_V = \mu^{-1}(0) = \left\{(E,F,v,x_1,x_2,\alpha)\in T^*X : x_1 + x_2 + v\otimes\alpha = 0\right\}.
    \]
    Like in the classical setting, one can obtain from each irreducible component of $Z_V$ a pair of SYTs $(T_1,T_2)$ by the sequences of Jordan types of $x_1,x_2$ on $E,F$, respectively. Let $\nu$ and $\nu'$ be the shapes of $T_1$ and $T_2$. We no longer have $x_1 + x_2 = 0$, so $\nu$ and $\nu'$ are not necessarily the same. However, they cannot be too different either.
    
    Since $x_1 = -x_2 - v\otimes\alpha$, one can actually show that $\CC[x_1] v = \CC[x_2]v =:U$. Then, on $V/U$, we do have $\Bar{x}_1 + \Bar{x}_2 = 0$, so $\Bar{x}_1$ and $\Bar{x}_2$ have the same Jordan type $\theta$. In addition, $U$ is a cyclic subspace generated by a vector $v$, so $\nu/\theta$ and $\nu'/\theta$ are horizontal strips. Let $\bfT_n$ be the set of such tripes $(\nu,\nu',\theta)$ where $|\nu| = |\nu'| = n$.

    Travkin \cite[Proposition 1]{travkin2009mirabolic} showed that these conditions are enough to classify irreducible components of $Z_V$. That is, we can index the irreducible components by triples $(T_1,T_2,\theta)$, where $T_1\in \SYT(\nu), T_2 \in \SYT(\nu')$, and $(\nu,\nu',\theta) \in \bfT_n$. Thus, we have a bijection between $RB_n$ and triples $(T_1,T_2,\theta)$ above. This is called mirabolic RSK.

    \begin{definition}[{\cite[Section 3.5]{travkin2009mirabolic}}]
        Let $\Tilde{w}\in RB_n$. The
        \textit{mirabolic Robinson--Schensted--Knuth correspondence} assigns to $\Tilde{w}$ a triple $(P(w),Q(w),\theta)$ as follows.
        
        Start with $P=Q=\varnothing$, and a separate row $r^@$ of $@ > n$. We have two stages.

        In Stage 1, for $i=1,\ldots,n$, insert $w(i)$ as follows.
        
        \begin{enumerate}
            \item If $i\in \beta$, insert $w(i)$ into $P$ similar to classical RS.
            \item If $i\notin \beta$, first insert $w(i)$ into $r^@$, then insert the replaced entry into $P$.
            \item Both cases create exactly one new box in $P$. Add a box in the same position in $Q$ and fill it with $i$.
        \end{enumerate}

        After Stage 1, we have our $Q(w) = Q$. Now we proceed to Stage 2.
        
        \begin{enumerate}
            \item[4.] Insert all entries in $r^@$ into $P$ from smallest to largest.
            \item[5.] $P(w)$ is $P$ restricted to the entries $1,\ldots,n$.
            \item[6.] $\theta$ is the shape of $P$ with the first row removed.
        \end{enumerate}
    \end{definition}

    \begin{example}
        Continuing Example \ref{ex:RB-orbit}, let us insert $(2431,\{1,4\})$. In Stage 1, we insert
        \[
        \varnothing~,~ \ctab{\begin{ytableau} 
            {@} & {@} & \cdots
        \end{ytableau}}
        ~\longrightarrow~
        \ctab{\begin{ytableau} 
            *(yellow) 2
        \end{ytableau}}
        ~,~ \ctab{\begin{ytableau} 
            {@} & {@} & \cdots
        \end{ytableau}}
        ~\longrightarrow~
        \ctab{\begin{ytableau} 
            2 & *(yellow) @
        \end{ytableau}}
        ~,~ \ctab{\begin{ytableau} 
            4 & {@} & \cdots
        \end{ytableau}}
        \]
        \[
        \longrightarrow~
        \ctab{\begin{ytableau} 
            2 & 4 \\
            *(yellow) @
        \end{ytableau}}
        ~,~ \ctab{\begin{ytableau} 
            3 & {@} & \cdots
        \end{ytableau}}
        ~\longrightarrow~
        \ctab{\begin{ytableau} 
            1 & 4 \\
            2 \\
            *(yellow) @
        \end{ytableau}}
        ~,~ \ctab{\begin{ytableau} 
            3 & {@} & \cdots
        \end{ytableau}}.
        \]
        Hence,
        \[
        Q(w) ~=~ \ctab{\begin{ytableau} 
            1 & 2 \\
            3 \\
            4
        \end{ytableau}}.
        \]
        Now for Stage 2, we insert $r^@$
        \[
        \ctab{\begin{ytableau} 
            1 & 4 \\
            2 \\
            @
        \end{ytableau}}
        ~,~ \ctab{\begin{ytableau} 
            3 & {@} & \cdots
        \end{ytableau}}
        ~\longrightarrow~
        \ctab{\begin{ytableau} 
            1 & 3 \\
            2 & 4 \\
            @
        \end{ytableau}}
        ~,~ \ctab{\begin{ytableau} 
            {@} & \cdots
        \end{ytableau}}
        ~\longrightarrow~
        \ctab{\begin{ytableau} 
            1 & 3 & @ & \cdots \\
            2 & 4 \\
            @
        \end{ytableau}}.
        \]
        Thus,
        \[
        P(w) ~=~ \ctab{\begin{ytableau} 
            1 & 3 \\
            2 & 4
        \end{ytableau}}, \qquad
        \theta ~=~\ctab{\begin{ytableau} 
            {} & {} \\
            {}
        \end{ytableau}}.
        \]
        Let us check the geometry. Recall that we can take
        \[
        E = e_{1}~)~e_{2}~)~e_3 ~)~e_4,\qquad
        F = e_{2}~)~ e_{4}~)~e_3~)~e_{1}, \qquad v= e_1 + e_2.
        \]
        Let
        \[
        x_1 = \left( \begin{matrix}
            0 & E_{12} & E_{13} & E_{14} \\
            0 & 0 & E_{23} & E_{24} \\
            0 & 0 & 0 & E_{34} \\
            0 & 0 & 0 & 0
        \end{matrix} \right) \in \bfrak_E, \qquad
        x_2 = \left( \begin{matrix}
            0 & 0 & 0 & 0 \\
            E'_{21} & 0 & E'_{23} & E'_{24} \\
            E'_{31} & 0 & 0 & 0 \\
            E'_{41} & 0 & E'_{43} & 0
        \end{matrix} \right) \in \bfrak_F.
        \]
        The condition $x_1+x_2+v\otimes\alpha = 0$ means
        \[
        \left( \begin{matrix}
            e^*_1 & E_{12} + e^*_2 & E_{13} + e^*_13 & E_{14} + e^*_4 \\
            E'_{21} + e^*_1 & e^*_2 & E_{23} + E'_{23} + e^*_3 & E_{24} + E'_{24} + e^*_4 \\
            E'_{31} & 0 & 0 & E_{34} \\
            E'_{41} & 0 & E'_{43} & 0
        \end{matrix} \right) = 0,
        \]
        which means
        \[
        e^*_1 = e^*_2 = E_{12} = E_{34} = E'_{21} = E'_{31} = E'_{41} = E'_{43} = 0.
        \]
        Thus, generic $x_1$ and $x_2$ have the form
        \[
        x_1 = \left( \begin{matrix}
            0 & 0 & E_{13} & E_{14} \\
            0 & 0 & E_{23} & E_{24} \\
            0 & 0 & 0 & 0 \\
            0 & 0 & 0 & 0
        \end{matrix} \right) \in \bfrak_E, \qquad
        x_2 = \left( \begin{matrix}
            0 & 0 & 0 & 0 \\
            0 & 0 & E'_{23} & E'_{24} \\
            0 & 0 & 0 & 0 \\
            0 & 0 & 0 & 0
        \end{matrix} \right) \in \bfrak_F.
        \]
        One can check that the sequences of Jordan types of $x_1$ and $x_2$ on $E$ and $F$ match $P$ and $Q$, respectively. In addition, $\CC[x_1]v = \CC[x_2]v = \CC(e_1 + e_2)$, and one can check that on $V/\CC(e_1 + e_2)$, the Jordan types of $x_1$ and $x_2$ are both $(2,1) = \theta$.
    \end{example}

\section{Finite orbits classification}

\subsection{Classification}

    Recall that irreducible polynomial representations are indexed by partitions. Let $D$ be the determinant representation, and recall that
    \[
    L(\lambda + k(1^n)) \simeq L(\lambda) \otimes D^{\otimes k}.
    \]
    In particular,
    \begin{align*}
         L(k^n)&=D^{\otimes k},\\
         L(k+1,k^{n-1})&=V\otimes D^{\otimes k},\\
         L((k+1)^{n-1},k)&=\bigwedge^{n-1}V\otimes D^{\otimes k}\simeq V^*\otimes D^{\otimes k+1}.
    \end{align*}
    Now we prove our classification theorem.

    \begin{thm}\label{thm:classify}
        Let $n\geq 2$ and let $L(\lambda)$ be an irreducible polynomial $\GL(V)$-representation. Then $\FlVV\times L(\lambda)$ has finitely many $G$-orbits if and only if $\lambda$ is one of the following
        \[
        \begin{array}{ll}
        (k^n),& k\geq1,\\
        (k+1,k^{n-1}),& k\geq0,\\
        ((k+1)^{n-1},k),& k\geq0.
        \end{array}
        \]
    \end{thm}

    \begin{proof}
        First, we claim that if $\FlVV\times W$ has finitely many $G$-orbits, then $\dim W\leq n$. Indeed, when $(E,F)$ are opposite flags, their common stabilizer is a maximal torus $T\simeq (\CC^\times)^n$. The $G$-orbits over the orbit of this pair of flags are in bijection with the $T$-orbits in $W$, which each has dimension at most $n$. Thus, if there are finitely many $G$-orbits, $W$ has to also have dimension at most $n$. By Weyl's dimension formula, $L(\lambda)$ has dimension at most $n$ only if $\lambda$ is one of $(k^n), (k+1,k^{n-1}), ((k+1)^{n-1},k)$ for some $k$.

        Now we check the candidates. It suffices to show that there are finitely many $T$-orbits. For $D^{\otimes k}$, $k\geq 1$, let $\delta \in D^{\otimes k}$. We have $(t_1,\ldots,t_n) \delta = ( t_1\cdots t_n)^k\delta$. The map $(t_1,\ldots,t_n)\mapsto(t_1\cdots t_n)^k$ is surjective on $\CC^\times$, so we have two $T$-orbits $\{0\} \cup \CC^\times$.

        For $V\otimes D^{\otimes k}$, we have $(t_1,\ldots,t_n) (e_i\otimes\delta) = t_i( t_1\cdots t_n)^k (e_i\otimes\delta)$. Thus, the torus can scale each coordinate independently, so the orbits depend only on the support. Hence, we have $2^n$ $T$-orbits. The argument for $V^*\otimes D^{\otimes k+1}$ is similar.

        Finally, for $\lambda = (0)$, $G$ acts trivially on $L(0)$, so there are infinitely many orbits.
    \end{proof}

\subsection{Orbit indices}

    For a pair of flags $(E,F)$ in relative position $w$, let
    \[
    A_w=\operatorname{span}_{\CC}\{E_{ij}:i\leq j \text{ and } w^{-1}(i)\leq w^{-1}(j)\},
    \]
    and $H_w = A_w^\times$. The $G$-orbits on $\FlVV\times W$ in which $(E,F)$ are in relative position $w$ are exactly the $H_w$-orbits on $W$. We will use this lemma.

   \begin{lemma}[{\cite[Lemma 1]{travkin2009mirabolic}}]\label{lem:orbit-submodule}
        Let $A \subset \operatorname{End}(V)$ be an associative algebra with identity, and let $A^\times$ be the multiplicative group of $A$. If the $A$-module $V$ has finitely many submodules, then $A^\times$-orbits in $V$ are in bijection with these submodules.
    \end{lemma}

    Recall that $RB_n$ is the set of pairs $\Tilde{w} = (w,\beta)$, where $w\in S_n$ and $\beta \subset [n]$ such that if $i\notin \beta$ and $j \in \beta$, then either $i > j$ or $w(i) > w(j)$. In other words, the marked squares are NW closed. 
    
    Let $\oRB_n$ be the set of opposite pairs $\Tilde{w} = (w,\beta)$, where $w\in S_n$ and $\beta \subset [n]$ such that if $i\notin \beta$ and $j \in \beta$, then either $i < j$ or $w(i) < w(j)$. In other words, the marked squares are SE closed.

    It is easy to see that $|RB_n| = |\oRB_n|$. Indeed, one can obtain $\Tilde{w}'\in \oRB_n$ from $\Tilde{w}\in RB_n$ by reflecting about the anti-diagonal.

    \begin{prop}
        We have the following bijections
        \begin{itemize}
            \item[(a)] $\GL(V)$-orbits on $\FlVV\times D^{\otimes k}$ $\longleftrightarrow$ pairs $(w,\eps)$ for $w\in S_n$, $\eps\in\{0,1\}$;
            \item[(b)] $\GL(V)$-orbits on $\FlVV\times (V\otimes D^{\otimes k})$ $\longleftrightarrow$ $RB_n$;
            \item[(c)] $\GL(V)$-orbits on $\FlVV\times (V^*\otimes D^{\otimes k+1})$ $\longleftrightarrow$ $\oRB_n$.
        \end{itemize}
    \end{prop}

    \begin{proof}
        For (a), $H_w$-submodules on $D^{\otimes k}$ are either $0$ or the whole $D^{\otimes k}$. We set $\eps = 0$ for $0$ and $\eps = 1$ for $D^{\otimes k}$.

        For (b), it is easy to see that $H_w$-submodules on $V$ correspond to RB permutations $(w,\beta)$. Adding a tensor factor of $D^{\otimes k}$ does not change this correspondence.

        The argument for (c) is the same, except that $\GL(V)$ acts on $V^*$ oppositely, so we use $\oRB_n$ instead of $RB_n$.
    \end{proof}

\section{Geometric RSK}

    Let us now denote
    \[
    X_k = \Fl(V)^2 \times D^{\otimes k},\quad
    X^+_k = \Fl(V)^2 \times (V\otimes D^{\otimes k}),\quad
    X^-_k = \Fl(V)^2 \times (V^*\otimes D^{\otimes k+1})
    \]
    with moment maps $\mu_k,\mu^+_k$, $\mu^-_k$, respectively, and let
    \[
    Z_k = \mu_k^{-1}(0),\quad Z^+_k = (\mu^+_k)^{-1}(0),\quad Z^-_k = (\mu^-_k)^{-1}(0)
    \]
    be the zero fibers.

\subsection{$D^{\otimes k}$}

    Let $\delta \in D^{\otimes k}$ and $A\in \gfrak$. Recall that $A\delta = \det(A)^k\delta$, so we have
    \[
    \exp(tA)\cdot \delta = (I+tA)\delta + t^2(\cdots) = \det(I+tA)^k\delta + t^2(\cdots) = \delta + kt\tr(A)\delta + t^2(\ldots).
    \]
    Thus,
    \[
    \left.\dfrac{d}{dt}\right|_{t=0} \exp(tA)\cdot \delta = k\tr(A)\delta,
    \]
    so $\mu_k$ is given by
    \[
    \langle \mu_k(E,F,\delta,x_1,x_2,\delta^*), A \rangle = \tr x_1 A + \tr x_2A + k\delta^*(\delta)\tr A = \langle x_1 + x_2 + k\delta^*(\delta)I ,A \rangle.
    \]
    Hence,
    \[
    \mu_k(E,F,\delta,x_1,x_2,\delta^*) = x_1 + x_2 + k\delta^*(\delta)I.
    \]
    Note that $x_1$ and $x_2$ both have zero on the diagonal, so the condition $x_1 + x_2 + k\delta^*(\delta)I = 0$ means that $x_1 + x_2 = 0$ and $\delta^*(\delta) = 0$. Thus,
    \[
    Z_k = \{(E,F,\delta,x,\delta^*)\in X_k\times \cN \times (D^{\otimes k})^*:\delta^*(\delta) = 0\}.
    \]
    Hence, we can describe $Z_k$-RSK from classical RSK as follows.

    \begin{thm}\label{thm:DkRSK}
        There is a correspondence on the irreducible components of $Z_k$
        \[
        (w,\eps)\longleftrightarrow (P(w),Q(w),\eps),
        \]
        where $w\longleftrightarrow (P(w),Q(w))$ is the classical RSK. In this correspondence, $P(w)$ and $Q(w)$ are the sequences of Jordan types of $x$ on $E$ and $F$, and $\eps = 0$ (resp. $\eps \neq 0$) if $\delta = 0$ (resp. $\delta \neq 0$).
    \end{thm}

\subsection{$V\otimes D^{\otimes k}$ and $V^* \otimes D^{\otimes k+1}$}

    Let
    \[
    Z_V = \{(E,F,v,x_1,x_2,\alpha):x_1+x_2+v\otimes\alpha = 0\}
    \]
    be the zero fiber in the setting of mirabolic RSK.

    \begin{thm}\label{thm:Z+-}
        We have $Z^+_k \simeq Z_V \simeq Z^-_k$.
    \end{thm}

    \begin{proof}
        Let $\delta \in D^{\otimes k}$ be a nonzero vector. Then, all elements in $V\otimes D^{\otimes k}$ can be written as $v\otimes \delta$ with the $G$-action by
        \[
        A\cdot (v\otimes \delta) = (\det A)^k (Av)\otimes \delta.
        \]
        Hence, by the same computation as above and in Section \ref{subsec:miraRSK}, $\mu^+_k$ is given by
        \[
        \mu^+_k(E,F,v\otimes \delta,x_1,x_2,\alpha\otimes \delta^*) = x_1 + x_2 + v\otimes\alpha + k\alpha(v) I.
        \]
        From the condition $ x_1 + x_2 + v\otimes\alpha + k\alpha(v) I = 0$, taking the trace gives $(1+kn)\alpha(v) = 0$, so $\alpha(v) = 0$. Hence, we have $x_1 + x_2 + v\otimes \alpha = 0$, so sending $(v\otimes \delta,\alpha\otimes \delta^*) \mapsto (v,\alpha)$ gives an isomorphism $Z^+_k \xrightarrow{\simeq} Z_V$.

        The argument for $V^*\otimes D^{\otimes k+1}$ is similar. The moment map $\mu^-_k$ is given by
        \[
        \mu^-_k(E,F,\alpha\otimes \delta,x_1,x_2,v\otimes \delta^*) = x_1 + x_2 - v\otimes\alpha + (k+1)\alpha(v) I.
        \]
        Thus, sending $(\alpha\otimes \delta^*,v\otimes \delta) \mapsto (v,-\alpha)$ gives an isomorphism $Z^-_k \xrightarrow{\simeq} Z_V$.
    \end{proof}

    \begin{cor}\label{cor:Z+-}
        \begin{itemize}
            \item There is a correspondence on the irreducible components of $Z^+_k$
            \[
            (w,\beta)\in RB_n \longleftrightarrow (T_1,T_2,\theta),
            \]
            given by mirabolic RSK. In this correspondence, $T_1$ and $T_2$ are the sequences of Jordan types of $x_1$ and $x_2$ on $E$ and $F$, respectively, and $\theta$ is the Jordan type of $x_1$ on $V/\CC[x_1]v$.
            \item There is a correspondence on the irreducible components of $Z^-_k$
            \[
            (w,\beta) \in \overline{RB}_n \longleftrightarrow (T_1,T_2,\theta'),
            \]
            where $T_1$ and $T_2$ are obtained from mirabolic RSK $(w,[n]\backslash \beta) \longleftrightarrow (T_1,T_2,\theta)$. Let $\nu$ and $\nu'$ be the shapes of $T_1$ and $T_2$, then
            \[
            \theta'_i = \min(\nu_i,\nu'_i) + \max(\nu_{i+1},\nu'_{i+1}) - \theta_i.
            \]
            In this correspondence, $T_1$ and $T_2$ are the sequences of Jordan types of $x_1$ and $x_2$ on $E$ and $F$, respectively, and $\theta'$ is the Jordan type of $x_1$ on $V/\CC[x_1]v$.
        \end{itemize}
    \end{cor}

\subsection{Fourier duality}

    Corollary \ref{cor:Z+-} is an instance of the Fourier exchange. For a vector space $W$, the canonical Fourier symplectic exchange is
    \[
    \cF_W:T^*W=W\times W^* \longrightarrow T^*(W^*)=W^*\times W,\qquad
    (w,\xi)\longmapsto(\xi,-w).
    \]
    This gives a Fourier bijection on the irreducible components of $T^*W$ and $T^*(W^*)$. Observe that $(V^*\otimes D^{\otimes k+1}) = (V\otimes D^{\otimes k})^* \otimes D^{\otimes 2k+1}$. Thus, after suitable rescaling, we have the following.

    \begin{thm}\label{thm:Fourier-V}
        The Fourier exchange $\cF$ gives an isomorphism
        \[
        \cF:Z_V^+\xrightarrow{\simeq} Z_V^-.
        \]
        On the $G$-orbits, $\cF$ sends
        \[
        (w,\beta) \in RB_n \mapsto (w,[n]\backslash \beta) \in \overline{RB}_n.
        \]
        On the triple $(T_1,T_2,\theta)$, $\cF$ sends
        \[
        (T_1,T_2,\theta) \mapsto (T_1,T_2,\theta'),
        \]
        where $\theta'_i = \min(\nu_i,\nu'_i) + \max(\nu_{i+1},\nu'_{i+1}) - \theta_i$, $\nu = \sh(T_1), \nu' = \sh(T_2)$.
    \end{thm}

    There is a second operation one can do, which is taking the dual of a flag. For a flag
    \[
    F = (0 = F_0 \subset F_1 \subset \cdots \subset F_n = V),
    \]
    its dual is
    \[
    F^\vee = (0 = F_n^\perp \subset F_{n-1}^\perp \subset \cdots \subset F_0^\perp = V^*).
    \]

    \begin{thm}\label{thm:Fourier-flag}
        The map
        \[
        \Delta:(E,F,v,x_1,x_2,\alpha) \mapsto (E^\vee,F^\vee,v,-x_1^*,-x_2^*,\alpha)
        \]
        extends to an isomorphism
        \[
        \Delta:Z^+_V \xrightarrow{\simeq} Z^-_{V^*}.
        \]
        On the $G$-orbits, $\Delta$ sends
        \[
        (w,\beta) \mapsto (w_0ww_0,w_0 \beta).
        \]
        On the triple $(T_1,T_2,\theta)$, $\Delta$ sends
        \[
        (T_1,T_2,\theta) \mapsto (\evac(T_1),\evac(T_2),\theta).
        \]
    \end{thm}

    We refer the reader to \cite{stanley2008promotion} for the definition of evacuation. Note that we have a commuting diagram
    \[\begin{tikzcd}[sep = scriptsize]
    	{Z^+_V} && {Z^-_V} \\
    	\\
    	{Z^-_{V^*}} && {Z^+_{V^*}}
    	\arrow["{\cF}", from=1-1, to=1-3]
    	\arrow["\Delta"', from=1-1, to=3-1]
    	\arrow["\Delta", from=1-3, to=3-3]
    	\arrow["{\cF}"', from=3-1, to=3-3]
        \arrow["{\mathcal{T}}"{description}, dotted, from=1-1, to=3-3]
    \end{tikzcd}\]
    in which the diagonal is Travkin's Fourier exchange on the whole $\FlVV \times V$ \cite[Section 5.3]{travkin2009mirabolic}.

    On the $G$-orbits, the diagram is
    \[\begin{tikzcd}[sep = scriptsize]
    	{(w,\beta)} && {(w,[n]\backslash \beta)} \\
    	\\
    	{(w_0ww_0,w_0\beta)} && {(w_0ww_0,[n]\backslash w_0\beta)}
    	\arrow["{\cF}", from=1-1, to=1-3]
    	\arrow["\Delta"', from=1-1, to=3-1]
    	\arrow["\Delta", from=1-3, to=3-3]
    	\arrow["{\cF}"', from=3-1, to=3-3]
        \arrow["{\mathcal{T}}"{description}, dotted, from=1-1, to=3-3]
    \end{tikzcd},\]
    and on the RSK triple, the diagram is
    \[\begin{tikzcd}[sep = scriptsize]
    	{(T_1,T_2,\theta)} && {(T_1,T_2,\theta')} \\
    	\\
    	{(\evac(T_1),\evac(T_2),\theta)} && {(\evac(T_1),\evac(T_2),\theta')}
    	\arrow["{\cF}", from=1-1, to=1-3]
    	\arrow["\Delta"', from=1-1, to=3-1]
    	\arrow["\Delta", from=1-3, to=3-3]
    	\arrow["{\cF}"', from=3-1, to=3-3]
        \arrow["{\mathcal{T}}"{description}, dotted, from=1-1, to=3-3]
    \end{tikzcd}.\]

\section{Mirabolic RSK as Up and Down operators}\label{sec:mira-UD}

    First, we review Viennot's shadow construction \cite{viennot2006forme}. Given some $\times$'s in an $n\times n$ grid. For this version, we place the source of light at the bottom right corner. Each $\times$ blocks the light to its NW. The boundary between blocked and unblocked area on the grid is the first shadow line. We remove the $\times$'s on this shadow line and repeat.

    \begin{example}
        Figure \ref{fig:shadow31862547} shows the shadow lines for $w = 31862547$.
        \begin{figure}[h!]
            \centering
            \includegraphics[scale = 1]{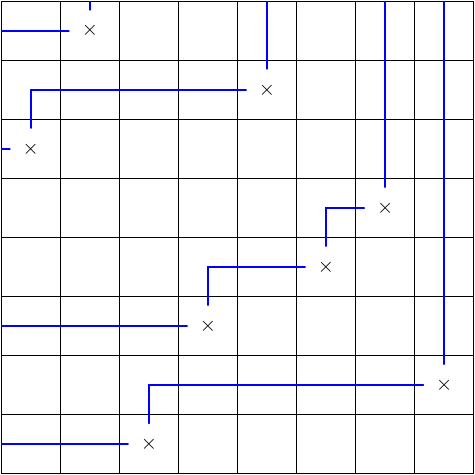}
            \caption{Shadow line}
            \label{fig:shadow31862547}
        \end{figure}
    \end{example}

    Let $P_n$ be the set of partial permutations of $[n]$. The following proposition gives a map $RB_n \longrightarrow P_n$. We leave the proof as an exercise to the reader. The idea of the inverse map is to place the light source at the opposite end of the grid.

    \begin{prop}
        There is a bijection $\Psi:RB_n \longrightarrow P_n$ as follows. Given an RB permutation $(w,\beta)$ drawn in an $n\times n$ grid. On the $\otimes$'s, draw the shadow lines and place a new $\times$ for each $\Big\ulcorner$-corner. The new $\times$'s and the old unmarked $\times$'s together form a partial permutation.

        For our incoming proof, let us also mark the two ends of each shadow line with a \textit{boundary} $\times$, place right above and left of the grid.
    \end{prop}

    \begin{example}\label{ex:Psi-map}
        Let $(w,\beta) = (53421,\{2,3,5\})$ shown in Figure \ref{fig:RBtoP}. The shadow lines have one new $\Big\ulcorner$ at $(1,3)$. Thus, the corresponding partial permutation is $1\rightarrow 5, 3\rightarrow 1, 4\rightarrow 2$. The boundary marks are colored green.
        \begin{figure}[h!]
            \centering
            \includegraphics[scale = 1]{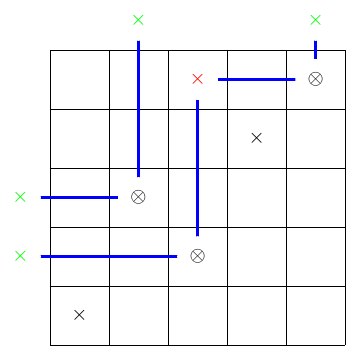}
            \caption{$\Psi$}
            \label{fig:RBtoP}
        \end{figure}
    \end{example}

    Next, we construct another map from the set of mirabolic triples $(T_1,T_2,\theta)$ to $P_n$ via the Up and Down operators.

    \begin{definition}
        Given a mirabolic triple $(T_1,T_2,\theta)$, where $T_1$ and $T_2$ are SYTs of shape $\nu$ and $\nu'$, for $|\nu| = |\nu'| = n$, $|\theta| = k$, and $\nu/\theta$ and $\nu'/\theta$ are horizontal strips. View $(T_1,T_2,\theta)$ as an oscillating tableau in $\OT(\varnothing,\varnothing,(1^n,-(n-k),n-k,(-1)^n))$. Define a map $\Xi : \OT(\varnothing,\varnothing,(1^n,-(n-k),n-k,(-1)^n)) \longrightarrow P_n$ as follows. Start with an $n\times n$ grid and add a boundary column at $j = -1$ and a boundary row at $i = -1$. Put $T_1$ along the East boundary from bottom up, $T_2$ along the North boundary from right to left, and $\theta$ at $(0,0)$. Complete the grid using Pak--Postnikov bijection $\Phi$ in Section \ref{subsec:PPUD} and put an $\times$ whenever $\Phi$ sends $T \in \OT(\lambda,\mu,(m,-n))$ into $\OT(\lambda,\mu,(-n+1,m-1))$. Here one of our $m$ and $n$ will be $1$. The $\times$'s inside the original grid is the desired partial permutation.
    \end{definition}

    \begin{example}\label{ex:Xi-map}
        Let
        \[
        T_1 = \ctab{\begin{ytableau}
            1 & 2 & 3 \\
            4 & 5
        \end{ytableau}}, \qquad T_2 = \ctab{\begin{ytableau}
            1 & 3 & 4 \\
            2 \\
            5
        \end{ytableau}}, \qquad \theta = \ctab{\begin{ytableau}
            {} & {} \\
            {}
        \end{ytableau}}.
        \]
        We can view
        \ytableausetup{boxsize=2.5em}
        \[
        (T_1,T_2,\theta) = \ctab{\begin{ytableau}
            1~\overline{12} & 2~\overline{10} & 3\overline{6}7\overline{9} \\
            4~\overline{11} & 5\overline{6} \\
            7\overline{8}
        \end{ytableau}}.
        \]
        The corresponding growth diagram is
        \begin{center}
            \includegraphics[scale = 0.7]{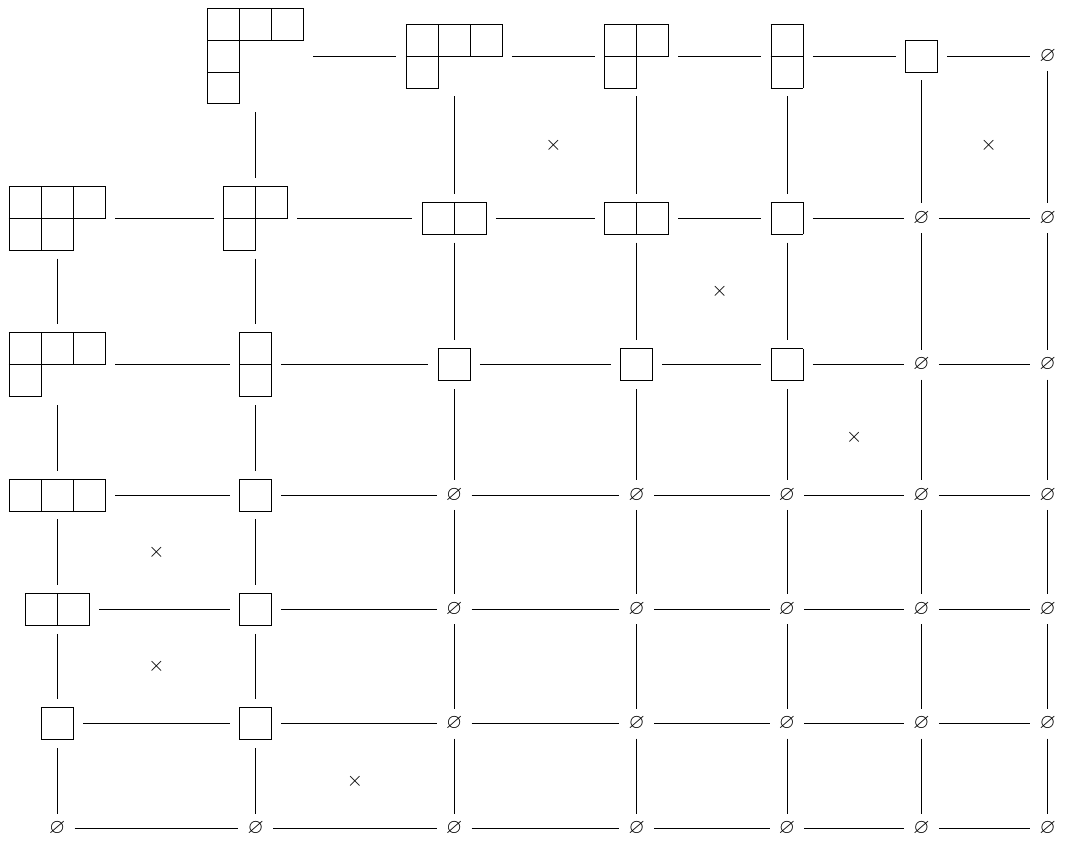}.
        \end{center}
        This gives us the partial permutation $1 \mapsto 5$, $3\mapsto 1$, $4\mapsto 2$.
    \end{example}
    \ytableausetup{boxsize=1.5em}

    The following is our main theorem of this section.

    \begin{thm}\label{thm:mira-UD}
        Let $\pi \in P_n$, let $(T_1,T_2,\theta)$ be the mirabolic RSK of $\Psi^{-1}(\pi)$, and $(T_1',T_2',\theta') = \Xi^{-1}(\pi)$. Then,
        \[
        T_1 = \evac(T_1'),\qquad T_2 = \evac(T_2'),\qquad \theta = \theta'.
        \]
    \end{thm}

    \begin{example}
        Mirabolic RSK sends $(w,\beta) = (53421,\{2,3,5\})$ in Example \ref{ex:Psi-map} to
        \[
        T_1 = \ctab{\begin{ytableau}
            1 & 2 & 5 \\
            3 & 4
        \end{ytableau}}, \qquad T_2 = \ctab{\begin{ytableau}
            1 & 3 & 4 \\
            2 \\
            5
        \end{ytableau}}, \qquad \theta = \ctab{\begin{ytableau}
            {} & {} \\
            {}
        \end{ytableau}}.
        \]
        One can check that $T_1$ and $T_2$ here are evacuations of the tableaux in Example \ref{ex:Xi-map}.
    \end{example}

    The following lemma, which follows easily from the toggling definition, is key to proving Theorem \ref{thm:mira-UD}.

    \begin{lemma}\label{lem:replace}
        Let $D_1$ be a $1\times n$ growth diagram in which the North boundary is a SYT of shape $\lambda$ arranged from right to left, and the East boundary is $\mu \subset \lambda$ such that $\lambda /\mu$ is a horizontal strip of size $k$
        \[\begin{tikzcd}[sep=scriptsize]
        	\lambda & {\lambda_{n-1}} & \cdots & {\lambda_1} & \varnothing \\
        	\mu & {} & \cdots & {} & \varnothing
        	\arrow[no head, from=1-1, to=1-2]
        	\arrow[no head, from=1-1, to=2-1]
        	\arrow[no head, from=1-2, to=1-3]
        	\arrow[no head, from=1-3, to=1-4]
        	\arrow[no head, from=1-4, to=1-5]
        	\arrow[no head, from=1-5, to=2-5]
        	\arrow[no head, from=2-1, to=2-2]
        	\arrow[no head, from=2-2, to=2-3]
        	\arrow[no head, from=2-3, to=2-4]
        	\arrow[no head, from=2-4, to=2-5]
        	\arrow[no head, from=1-2, to=2-2]
        	\arrow[no head, from=1-4, to=2-4]
        \end{tikzcd}.\]
        Let $D_2$ be a $k\times n$ growth diagram with the North boundary is the same as $D_1$, and the East boundary is a sequence $\lambda \supset \nu_1 \supset \cdots \supset \nu_{k-1} \supset \mu$ such that in each step, we remove the rightmost box possible
        \[\begin{tikzcd}[sep=scriptsize]
        	\lambda & {\lambda_{n-1}} & \cdots & {\lambda_1} & \varnothing \\
        	{\nu_1} & {} & \cdots & {} & \varnothing \\
        	\vdots & \vdots && \vdots & \vdots \\
        	{\nu_{k-1}} & {} & \cdots & {} & \varnothing \\
        	\mu & {} & \cdots & {} & \varnothing
        	\arrow[no head, from=1-1, to=1-2]
        	\arrow[no head, from=1-1, to=2-1]
        	\arrow[no head, from=1-2, to=1-3]
        	\arrow[no head, from=1-2, to=2-2]
        	\arrow[no head, from=1-3, to=1-4]
        	\arrow[no head, from=1-4, to=1-5]
        	\arrow[no head, from=1-4, to=2-4]
        	\arrow[no head, from=1-5, to=2-5]
        	\arrow[no head, from=2-1, to=2-2]
        	\arrow[no head, from=2-1, to=3-1]
        	\arrow[no head, from=2-2, to=2-3]
        	\arrow[no head, from=2-2, to=3-2]
        	\arrow[no head, from=2-3, to=2-4]
        	\arrow[no head, from=2-4, to=2-5]
        	\arrow[no head, from=2-4, to=3-4]
        	\arrow[no head, from=2-5, to=3-5]
        	\arrow[no head, from=3-1, to=4-1]
        	\arrow[no head, from=3-2, to=4-2]
        	\arrow[no head, from=3-4, to=4-4]
        	\arrow[no head, from=3-5, to=4-5]
        	\arrow[no head, from=4-1, to=4-2]
        	\arrow[no head, from=4-1, to=5-1]
        	\arrow[no head, from=4-2, to=4-3]
        	\arrow[no head, from=4-2, to=5-2]
        	\arrow[no head, from=4-3, to=4-4]
        	\arrow[no head, from=4-4, to=4-5]
        	\arrow[no head, from=4-4, to=5-4]
        	\arrow[no head, from=4-5, to=5-5]
        	\arrow[no head, from=5-1, to=5-2]
        	\arrow[no head, from=5-2, to=5-3]
        	\arrow[no head, from=5-3, to=5-4]
        	\arrow[no head, from=5-4, to=5-5]
        \end{tikzcd}.\]
        Then, $D_1$ and $D_2$ have the same South boundary. Furthermore, if $D_1$ have $\times$'s in columns $c_1,\ldots,c_k$, then $D_2$ have $\times$'s in exactly $(1,c_1),\ldots,(k,c_k)$.

        The transpose statement is also true.
    \end{lemma}

    \begin{example}
        Let $D_1$ be the top row in Example \ref{ex:Xi-map}
        \begin{center}
            \includegraphics[scale = 0.7]{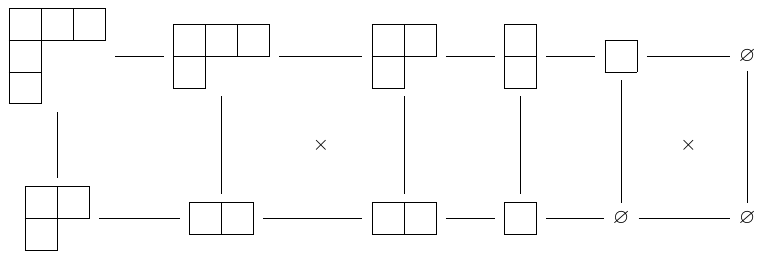}.
        \end{center}
        Then, $D_2$ is the following diagram
        \begin{center}
            \includegraphics[scale = 0.7]{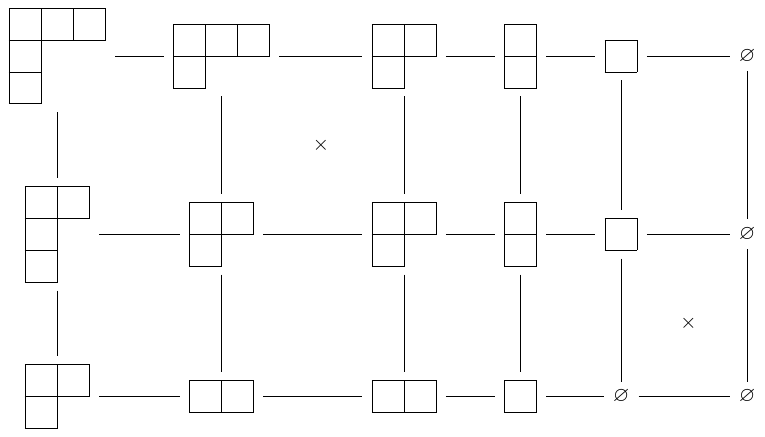}.
        \end{center}
    \end{example}

    \begin{proof}[Proof of Theorem \ref{thm:mira-UD}]
        Replacing the top row and leftmost column in $\Xi$ by the new $k\times n$ row and $n\times k$ column in Lemma \ref{lem:replace}, we have a new permutation. This is exactly the permutation $w_1$ in \cite[Lemma 3]{travkin2009mirabolic} for $(w_0ww_0,[n]\backslash w_0\beta)$. It follows from the toggling on $\theta$ in Theorem \ref{thm:Fourier-V} that $(T_1',T_2',\theta')$ is the mirabolic triple for $Z^-_{V^*}$. Theorem \ref{thm:mira-UD} then follows from Theorem \ref{thm:Fourier-flag}.
    \end{proof}

\section{Homology of mirabolic Steinberg variety}

    Recall from Theorem \ref{thm:PPrep} that, for each $0\leq k\leq n$, there is a natural Pak--Postnikov representation 
    \[
    M(n,(-k,k),n) = \bigoplus_{\substack{|\nu| = |\nu'| = n,|\theta| = n-k \\ \nu/\theta,\nu'/\theta\text{ horizontal strips}}} \pi_\nu \boxtimes \pi^*_{\nu'}.
    \]
    We now want to find these representations inside $H(Z_V)$. We have
    \[
    H(Z_V) = \bigoplus_w \CC\cdot[Z_{(T_1,T_2,\theta)}],
    \]
    Let
    \[
    Z_{\leq k}=\{(E,F,v;x_1,x_2,\alpha)\in Z_V:\dim \CC[x_1]v\leq k\},\qquad Z_{\le-1}=\varnothing,
    \]
    and define the filtration
    \[
    \scrF_k H(Z_V)=\im \left(H(Z_{\leq k}) \longrightarrow H(Z_V)\right), \qquad -1\leq k\leq n.
    \]
    \begin{lemma}
        Each $Z_{\leq k}$ is closed, and
        \[
        \scrF_k H(Z_V)=
        \bigoplus_{(T_1,T_2,\theta):n - |\theta| \leq k}
                     \CC\cdot[Z_{(T_1,T_2,\theta)}].
        \]
        Consequently $\gr_kH(Z_V)=
        \scrF_k H(Z_V)/\scrF_{k-1}H(Z_V)$ has a component basis indexed by the triples with $|\theta|=n-r$.
    \end{lemma}

    \begin{proof}
        The condition $\dim\CC[x_1]v \leq k$ is the same as saying the rank of the matrix
        \[
        \begin{bmatrix}
            v & x_1v & \cdots & x_1^{n-1}v
        \end{bmatrix}
        \]
        is at most $k$. Thus, $Z_{\leq k}$ is closed, and hence the entire component $Z_{(T_1,T_2,\theta)}$ lies in $Z_{\leq k}$ if and only if $n - |\theta| \leq k$.
    \end{proof}

    Our main theorem is the following.

    \begin{thm}\label{thm:gr-PP}
        For every $0\leq k\leq n$, there is an isomorphism of $S_n\times S_n$-representations
        \[
        \gr_k H(Z_V) \simeq M(n,(-k,k),n) \simeq \bigoplus_{(\nu,\nu',\theta)\in \bfT_n, |\theta| = n-k} \pi_\nu \boxtimes \pi^*_{\nu'}.
        \]
    \end{thm}

    First, we construct the $S_n\times S_n$-action. Let
    \[
    Y_V=\{(x_1,x_2,v,\alpha)\in\cN\times \cN\times V\times V^*: x_1+x_2+v\otimes\alpha = 0\},
    \]
    and $\Tilde{\cN} = T^*\Fl(V)$ with the Springer resolution $s:\Tilde{\cN} \rightarrow \cN$. Then we have a commuting diagram of forgetful maps
    \[\begin{tikzcd}
    	{Z_V} & {\Tilde{\mathcal{N}} \times \Tilde{\mathcal{N}}} \\
    	{Y_V} & {\mathcal{N} \times \mathcal{N}}
    	\arrow["g", from=1-1, to=1-2]
    	\arrow["p"', from=1-1, to=2-1]
    	\arrow["s\times s", from=1-2, to=2-2]
    	\arrow["f"', from=2-1, to=2-2]
    \end{tikzcd}.\]
    The map $p$ is proper and $p^{-1}(x_1,x_2,v,\beta)=\cB_{x_1}\times\cB_{x_2}$. Recall that if $x$ has Jordan type $\lambda$, then $H(\cB_x) \simeq \pi_\lambda$.

    \begin{prop}\label{prop:pushforward-action}
        There are two commuting $S_n$-actions on $H(Z_V)$. For every closed $A\subseteq Y_V$, the pushforward
        \[
        H(p^{-1}A)\longrightarrow H(Z_V)
        \]
        is equivariant. In particular, $\scrF$ is a filtration by $S_n\times S_n$-subrepresentations.
    \end{prop}

    \begin{proof}
        Recall the classical Steinberg variety $Z = \Tilde{\cN} \times_{\cN} \Tilde{\cN}$ and the convolution product in Section \ref{subsec:convo-hom-st}. In particular, the convolution product gives $H_*(Z)$ an associative algebra structure, and $H(Z) \simeq \CC[S_n]$. Let $[\Delta]$ be the unit of this subalgebra. This corresponds to the pairs of flags with relative position $w = \operatorname{id}$.

        Let $\Tilde{Y}_V = \gfrak\times\gfrak\times V\times V^* \supset Y_V$. Let $Z_A = p^{-1}(A)$ and $j_A:Z_A\hookrightarrow Z_V$ be the inclusion map. To set up the convolution product, let
        \[
        S_A:=\{(x_0,x_1,x_2,x_3,y) \in \Tilde{\cN}^4 \times \Tilde{Y}_V:\ (x_0,x_1)\in Z, (x_1,x_2,y)\in Z_A, (x_2,x_3)\in Z\},
        \]
        and $q:(x_0,x_1,x_2,x_3,y) \mapsto (x_0,x_3,y)$ be the projection map. We have $q(S_A) \subseteq Z_A$ since $s(x_0)=s(x_1), s(x_2)=s(x_3)$. Moreover, we have
        \[
        S_A\simeq Z_A\times_{\cN^2}(\Tilde{\cN}\times \Tilde{\cN}),
        \]
        so $q$ is the base change of the proper map $s\times s$, and hence is itself proper.
        
        Define the projections
        \[
        \pi_{01}: \Tilde{\cN}^4 \times \Tilde{Y}_V\rightarrow \Tilde{\cN}^2,\qquad
        \pi_{12Y}: \Tilde{\cN}^4 \times \Tilde{Y}_V\rightarrow \Tilde{\cN}^2 \times \Tilde{Y}_V,\qquad
        \pi_{23}: \Tilde{\cN}^4 \times \Tilde{Y}_V\rightarrow \Tilde{\cN}^2.
        \]
        For $x,y\in H(Z)$ and $z\in H_j(Z_A)$, let
        \[
        \widetilde x=\pi_{01}^*x,\qquad
        \widetilde z=\pi_{12Y}^*z,\qquad
        \widetilde y=\pi_{23}^*y,
        \]
        where the pullbacks are exterior products with the fundamental classes of the missing smooth factors. Their supports intersect in $S_A$, so we can define
        \[
        x\ast z\ast y = q_*(\Tilde{x} \cap \Tilde{z} \cap \Tilde{y}) \in H(Z_A).
        \]
        It is not hard to check that this preserves the degree of $z$. Thus, we have a left and a right action of $H(Z)$ on $H(Z_A)$ by
        \[
        x\cdot z = x\ast z\ast [\Delta],\qquad z\cdot y = [\Delta] \ast z\ast y.
        \]
        It is not hard to check that these left and right actions make $H(Z_A)$ an $H(Z)$-bimodule. Now we check that the pushforward is equivariant. By \cite[Section 2.6]{chriss1997representation}, 
        \[
        \Tilde{x} \cap \widetilde{(j_A)_*z} \cap \Tilde{y} = \widetilde{j_A}_*(\Tilde{x} \cap \Tilde{z} \cap \Tilde{y}),
        \]
        where $\widetilde{j_A}:S_A \hookrightarrow S_{Z_V}$ is the closed embedding. Since $q|_{S_{Z_V}}\circ\widetilde{j_A}=j_A\circ q|_{S_{A}}$,
        \[
        x\ast (j_A)_*z \ast y = (j_A)_*(x\ast z\ast y).
        \]
        Thus, $H(p^{-1}A)\to H(Z_V)$ is $S_n\times S_n$-equivariant.

        Finally, let $Y_{\leq k} = \{(x_1,x_2,v,\alpha)\in Y_V: \dim \CC[x_1]v \leq k\} \subset Y_V$. Then $Y_{\leq k}$ is closed, and $Z_{\leq k} = p^{-1}(Y_{\leq k})$. Then, $\scrF_k H(Z_V)=\im \left(H(Z_{\leq k}) \longrightarrow H(Z_V)\right)$ is the image of an equivariant map and hence is an $S_n\times S_n$-subrepresentation. 
    \end{proof}

    For a triple $\tau=(\nu,\nu',\theta)\in \bfT_n$, let
    \[
    Y_\tau=\{(x_1,x_2,v,\alpha)\in Y_V:J(x_1) = \nu,J(x_2) = \nu', J(x_1|_{V/\CC[x_1]v}) = \theta\}.
    \]
    Let $d = \dim Z_V$ and $b_\tau = b(\nu) + b(\nu')$, where $b(\nu)$ and $b(\nu')$ are the dimensions of $\cB_{\nu}$ and $\cB_{\nu'}$. For $e=(x_1,x_2,v,\alpha)\in Y_\tau$, we have $F_e=p^{-1}(e)=\cB_{x_1}\times\cB_{x_2}$, which is equidimensional of dimension $b_\tau$. Thus, we have $\dim Y_\tau\leq d-b_\tau$.

    \begin{lemma}\label{lem:trivial}
        Let $U\subseteq Y_\tau$ be smooth, connected and locally closed, with $\dim U=d-b_\tau$, and $p\colon Z_U=p^{-1}(U)\to U$ is a locally trivial fibration. Then for every $e\in U$,
        restriction to the fiber gives an $S_n\times S_n$-equivariant isomorphism
        \[
        H_{2d}(Z_U)\xrightarrow{\simeq}H_{2b_\tau}(F_e)\simeq\pi_\nu\boxtimes\pi_{\nu'}^*.
        \]
    \end{lemma}
    
    \begin{proof}
        Let $\mathcal L$ be the local system on $U$ with stalks $H_{2b_\tau}(F_{e'})$ for $e'\in U$. Since $Z_U\to U$ is a locally trivial fibration with fibers of
        dimension $b_\tau$ over a connected base of dimension $d-b_\tau$, we have $H_{2d}(Z_U)= \Gamma(U,\mathcal L)=H_{2b_\tau}(F_e)^{\pi_1(U,e)}$.
        
        Since $H(Z)$ acts on $H(Z_{U'})$ for every open $U'\subseteq U$, compatibly with restriction, it acts on $\cL$ by endomorphisms of local systems, so the monodromy commutes with the $H(Z)$-action. On the stalk at $e$, convolution with the fixed $Y_V$ is the Springer action on $H(\cB_{x_1})\boxtimes H(\cB_{x_2})$. This gives $H_{2b_\tau}(F_e) \simeq \pi_\nu\boxtimes\pi_{\nu'}^*$, which is
        irreducible.
        
        Finally, by Schur's lemma, each $\gamma\in\pi_1(U,e)$ acts by a scalar $c(\gamma)$. On the other hand, $\dim Z_U=d$, so $Z_U$ has a $d$-dimensional component, and its fundamental class is a nonzero invariant. Therefore $c(\gamma)=1$ for all $\gamma$, so $H_{2d}(Z_U) \simeq H_{2b_\tau}(F_e)$.
    \end{proof}

    \begin{lemma}\label{lem:irreducibility}
        $Y_\tau$ has a unique irreducible component $Y_\tau^{\circ}$ of dimension $d-b_\tau$, and all
        its other components have smaller dimension. Moreover, every irreducible component $C$ of
        $Z_V$ whose generic point has Jordan types $\tau$ satisfies $\overline{p(C)} = \overline{Y_\tau^\circ}$.
    \end{lemma}
    
    \begin{proof}
        Let $C$ be a component of $Z_V$ with generic Jordan types $\tau$. Then $p(C)\cap Y_\tau$ is
        dense in $\overline{p(C)}$, and its fibers have dimension at most $b_\tau$, so $\dim \overline{p(C)} \geq d-b_\tau \geq \dim Y_\tau$. Hence, $p(C)\cap Y_\tau$ is dense in a component of $Y_\tau$ of dimension exactly $d-b_\tau$.

       Let $C_1,\dots,C_r$ be the components of $Y_\tau$ of dimension $d-b_\tau$. We will show that $r = 1$. For each $i$, choose $U_i\subseteq C_i\setminus\bigcup_{j\ne i}C_j$ open, dense, smooth and connected with $p$ locally trivial over $U_i$. By Lemma~\ref{lem:trivial}, $Z_{U_i}$ is pure of dimension $d$ and has $|\SYT(\nu)|\cdot|\SYT(\nu')|$ components. Their closures are components of $Z_V$ with generic Jordan types $\tau$. They are also pairwise distinct, so together they contribute $r\cdot |\SYT(\nu)|\cdot |\SYT(\nu')|$ components in $Z_V$ with generic Jordan types $\tau$. By \cite[Proposition 1]{travkin2009mirabolic}, we have $r = 1$.
    \end{proof}
    
    Let $S_\tau=\overline{Y^\circ_\tau}$, and
    \[
    Z_{\leq \tau} = p^{-1}\left(\bigcup_{S_\sigma\subseteq S_\tau}S_\sigma\right),\qquad Z_{< \tau} = p^{-1}\left(\bigcup_{S_\sigma\subsetneq S_\tau}S_\sigma\right).
    \]
    We further refine $\cF$ by
    \[
    \cG_{\leq \tau} = \im\left(H(Z_{\leq \tau}) \rightarrow H(Z_V)\right), \qquad \cG_{< \tau} = \im\left(H(Z_{< \tau}) \rightarrow H(Z_V)\right).
    \]
    \begin{lemma}\label{lem:tau-filtration}
        For $\tau=(\nu,\nu',\theta)$ there is an isomorphism
        \[
        \cG_{\leq\tau}/\cG_{<\tau} \simeq \pi_\nu\boxtimes (\pi_{\nu'})^*.
        \]
    \end{lemma}

    \begin{proof}
        Choose $U_\tau\subseteq Y_\tau^\circ$ open, dense, smooth and connected such that $U_\tau$ avoids the other components of $S_\tau$, and $p$ is locally trivial over $U_\tau$. Then $U_\tau$ is open in $S_\tau$. Pushforward along closed embeddings is injective in top degree, so $G_{\leq\tau}$ has a basis of the classes $[C]$ of the components $C$ of $Z_V$ with $p(C)\subseteq S_\tau$. By Lemma \ref{lem:irreducibility}, these are exactly the components with $S_{\sigma(C)}\subseteq S_\tau$, where $\sigma(C)$ is the generic Jordan types of $C$. Similarly, $G_{<\tau}$ is spanned by the $[C]$ with $S_{\sigma(C)}\subsetneq S_\tau$. Hence, we have an isomorphism
        \[
        \cG_{\leq\tau} /\cG_{<\tau}\xrightarrow{\simeq} H(Z_{U_\tau}),
        \]
        which by Lemma \ref{lem:trivial} is $\pi_\nu\boxtimes\pi_{\nu'}^*$.
    \end{proof}

    \begin{proof}[Proof of Theorem \ref{thm:gr-PP}]
        Since $Z_{\le k}$ is closed, $S_\sigma\subseteq S_\tau$ implies $n-|\theta_\sigma|\le n-|\theta_\tau|$. Order the elements in $\bfT_n = \{\tau_1,\tau_2,\ldots\}$ first by $n-|\theta|$ and then by any order that respects inclusion of $S_\tau$'s. Let $\cH_{\leq i}=\sum_{j\leq i}\cG_{\leq \tau_j}$. By Lemma \ref{lem:irreducibility}, it is spanned by $\{[C]:\sigma(C) \in \{\tau_1,\dots,\tau_i\}\}$. Thus, $\cH$ is a filtration by subrepresentations that refines $\cF$, and $\cH_{\leq i}/\cH_{<i}\simeq\cG_{\leq \tau_i}/\cG_{<\tau_i}$. By Lemma~\ref{lem:tau-filtration},
        \[
        \gr_k H(Z_V) \simeq \bigoplus_{(\nu,\nu',\theta)\in \bfT_n, |\theta| = n-k} \pi_\nu \boxtimes \pi^*_{\nu'} \simeq M(n,(-k,k),n).
        \]
    \end{proof}

\section*{AI declaration}

    All ideas and content of this paper was formulated by the author. GPT-6 Astra and Claude Opus 5.5 were used to review literature, do RSK experiments, and fact-check technical lemmas.

\bibliography{bibliography}
\bibliographystyle{alpha}

\end{document}